\documentclass[12pt,final]{elsarticle}
\usepackage{amsmath,amsthm,amssymb}
\usepackage{a4wide}
\usepackage[latin1]{inputenc}
\usepackage{latexsym}
\usepackage[mathscr]{eucal}
\usepackage{epsfig}
\usepackage{amsfonts}
\usepackage{graphicx}
\usepackage{xcolor}

\newcommand\disk{\mathaccent"7017 }
\newcommand{\NN}{\mathbb{N}}
\newcommand{\ZZ}{\mathbb{Z}}
\newcommand{\RR}{\mathbb{R}}
\newcommand{\CC}{\mathbb{C}}
\newcommand{\TT}{\mathbb{T}}
\newcommand{\PP}{\mathbb{P}}
\newcommand{\LL}{\mathbb{L}}
\newcommand{\Ll}{\mathcal{L}}
\newcommand{\DD}{\mathbb{D}}
\newcommand{\EE}{\mathbb{E}}

\newcommand{\QQ}{\mathbb{Q}}

\newcommand{\Span}{\mathrm{span}\,}
\renewcommand{\v}[1]{\mathbf{#1}}

\newcommand{\arc}{\mathrm{arc}}

\newcommand{\supp}{\mathrm{supp}}
\newtheorem{theorem}{Theorem}[section]
\newtheorem{lemma}[theorem]{Lemma}
\newtheorem{corollary}[theorem]{Corollary}
\newtheorem{definition}[theorem]{Definition}
\newtheorem{remark}[theorem]{Remark}

\begin{document}
\begin{frontmatter}


\title{Zeros of quasi-paraorthogonal polynomials and positive quadrature}

\author[LEUVEN]{Adhemar Bultheel\corref{cor1}}

\address[LEUVEN]{Department of Computer Science, KU Leuven, Celestijnenlaan 200 A, B-3001
Heverlee, Belgium.}

\ead{adhemar.bultheel@cs.kuleuven.be}

\cortext[cor1]{Corresponding author.}

\author[ULL]{Ruym\'an Cruz-Barroso\fnref{label2}}

\ead{rcruzb@ull.es}

\author[ULL]{Carlos D\'{\i}az Mendoza\fnref{label2}}

\ead{cjdiaz@ull.es}

\fntext[label2]{This work is partially supported by Ministerio de Econom\'{\i}a y
Competitividad of Spain under grant MTM 2015-71352-P.}

\address[ULL]{Department of Mathematical Analysis, La Laguna
University, 38271 La Laguna, Tenerife. Canary Islands. Spain.\\[3mm]
{Dedicated to the memory of Luc Wuytack.}
}

\begin{abstract}
In this paper we illustrate that paraorthogonality on the unit circle $\TT$ is the
counterpart to orthogonality on $\RR$ when we are interested in the spectral properties.
We characterize quasi-paraorthogonal polynomials on the unit circle as the analogues of the
quasi-orthogonal polynomials on $\RR$. We analyze the possibilities of preselecting some
of its zeros, in order to build positive quadrature formulas with prefixed nodes and
maximal domain of validity.
These quadrature formulas on the unit circle are illustrated numerically.
\end{abstract}


\end{frontmatter}

\section{Introduction}\label{secIntro}
It is very well known that the orthogonal polynomials $\{\rho_n\}_{n\geq 0}$ associated with a
positive Borel measure $\mu$ supported on the real line (OPRL) have many properties that allow us,
both to emulate the measure, and describe its support.
For example, the zeros of $\rho_n$ are simple, they lie on the open convex hull of
$\supp(\mu)$ and there is one zero on the closure of each connected component of the
complement of $\supp(\mu)$, at most.
Moreover, the zeros of $\rho_{n+1}$ and $\rho_n$ interlace and in general,
they end up filling the whole support.
These zeros describe $\mu$ as $n \rightarrow \infty$ because they allow us to construct
a sequence of measures $\{\mu_n\}_{n\geq 1}$ that converges weakly to $\mu$.

Orthogonal polynomials on the unit circle $\TT=\{z \in \CC : |z|=1 \}$ (OPUC),
or Szeg\H{o} polynomials, were introduced in \cite[Chapter 11]{Sz}.
Since then, they have been widely studied, not only because of their own interest
\cite{SimonBk}, but also in many applications
such as the trigonometric moment problem \cite{Akh}, complex approximation \cite{Walsh},
probability and statistics \cite{GS},
prediction theory \cite{WM},
systems theory, networks, circuits and scattering \cite{DD},
signal processing \cite{Del}, and many more,
but also because of their intimate connection with OPRL (see e.g. \cite{Bul,CDPPablo}).
However, these polynomials present important differences with respect to OPRL,
in particular, concerning the above properties, since their zeros are located
outside of the support of the measure.
As a consequence, these zeros cannot be used for interpolation processes on the unit
circle, as nodes for quadrature formulas on the unit circle, and they lack
importance in the resolution of the trigonometric moment problem, among other
relevant properties that hold for OPRL.
These drawbacks were not solved until the concept of paraorthogonality (\cite{JNT},
see also \cite[Theorem III]{Ger}) was introduced in this setting.
These paraorthogonal polynomials (POPUC) have their zeros in the support of the
measure when the support is $\TT$. In this context, focussing on their spectral
properties, quadrature, and associated problems, these polynomials are a natural
counterpart on the unit circle of the orthogonal polynomials on the real line.
Some of the properties that we have mentioned above and many others formulated
for a measure supported on the real line, have been adapted for measures on the
unit circle, see e.g. \cite{CMV02,G,S,Wong}.

Since the spectra of the orthogonal and paraorthogonal polynomials are directly
connected to their respective type of orthogonality,
we shall also look at what can be said about the zeros as we relax some of the
orthogonality conditions
of these polynomials, leading to quasi-orthogonal polynomials on the real line
and the quasi-paraorthogonal counterpart on the unit circle.

In the case of the real line, the real polynomials $Q_{n,d}$ of degree $n$
that are orthogonal to $\PP_{n-d-1}$ (the vector space of all polynomials of
degree at most $n-d-1$),
$0 \leq d < n$ are well studied and they are called quasi-orthogonal polynomials
(QOPRL) of degree $n$ and order $d$ or $(n,d)$-QOPRL.
Clearly the orthogonal polynomial $\rho_n$ appears as a special case for $d=0$: $\rho_n=Q_{n,0}$.
The QOPRL $Q_{n,d}$ have at least $n-d$ changes of sign on the open convex hull of $\supp(\mu)$.
This concept was first introduced in 1923 by M. Riesz (\cite{Rie}) for $d=1$ in
relation to moment problems, and then by L. Fej\`er ($d=2$) and by J.A. Shohat
($d\geq 1$) in 1933 and 1937, respectively (\cite{Fej,Sho}), in the context of
quadrature formulas (see also \cite{Chi2,Dic,Dra}).

The unit circle analogue of these QOPRL should naturally be called quasi-paraorthogonal
polynomials (QPOPUC). Their properties however, are less developed in the literature.
This is one of the main purposes of this paper: to define and analyse the spectral
properties of these QPOPUC and how these relate to their orthogonality conditions.

The outline of the paper is as follows.
In Section~\ref{secQPOPUC} we introduce the concept of quasi-paraorthogonal polynomials on the unit circle, where  we shall generalize the classical concept of paraorthogonality to QPOPUC of order $2\ell+1$ by requiring
orthogonality to spaces of decreasing dimension. The classical concept of paraorthogonality corresponds to quasi-paraorthgonality of order 1.
We prove the main result in Section~\ref{secZero}: a link between the number of
paraorthogonality conditions and the number of zeros of odd multiplicity in the
interior of the support of the measure, when it is supported on an arc of the unit circle.
We also investigate how we can use the free parameters to prefix some of the zeros.
As an application we consider in Section~\ref{secQF} how this can be used to construct
positive Szeg\H{o}-type
quadrature formulas on an arc or on $\TT$ that are the analogues of the Gauss-type
formulas on an interval or $\RR$.
Finally, in Section~\ref{secNum} we shall include some numerical experiments to illustrate
 these quadrature formulas.

We finish this introductory section by defining some notation and abbreviations that
are used throughout this paper.
We denote by $\PP$ the vector space of all polynomials and by $\PP_n=\Span\{z^k:k=0,\ldots,n\}$
the $(n+1)$-dimensional space of polynomials of degree at most $n$.
The vector space of all Laurent polynomials is $\LL=\Span\left\{ z^k : k \in \ZZ \right\}$,
with subspaces
$\LL_{p,q}=\Span\{z^{p},\ldots,z^q\}$, $p\le q$, $p,q\in\ZZ$.
$\RR$ is the real line and $\TT$ the complex unit circle and $\mu$ is a
positive Borel measure that is supported on a subset $S$ of $\RR$ or $\TT$.
If $S \subsetneq \RR$ or $S \subsetneq \TT$ we assume it is compact and connected, i.e.\ an interval on $\RR$ or an arc on $\TT$.  We denote by $\langle \cdot, \cdot \rangle$ the inner product induced by $\mu$
and we assume that $\PP$  belongs to the corresponding Hilbert space $L_\mu^2$.
Orthogonality will always refer to this space.
The acronyms OPRL and OPUC are used for orthogonal polynomials on the real line and the unit circle respectively.
POPUC refers to paraorthogonal polynomials on the unit circle and if the letter Q is
added in the beginning, it refers to quasi-versions, thus QOPRL means quasi-orthogonal
polynomials on the real line and QPOPUC denotes
quasi-paraorthogonal polynomials on the unit circle. The precise definitions follow below.
By $\DD$ and $\EE$ we denote the interior and exterior of the closed unit disk,
respectively and by $\lfloor x\rfloor$ the integer part of $x\in \RR$.
$P^*(z)= z^n\overline{P\left(1/\overline{z}\right)}$  is the reciprocal (or reversed) polynomial of $P\in\PP_n$.
Throughout the paper, $\rho_n$ will denote the monic orthogonal polynomial of degree $n$ in $L^2_\mu$.
For $S\subseteq\TT$, the Schur parameters for the OPUC
will be denoted by $\{\delta_k\}_{k=0}^\infty$ where $\delta_0=1$ and
$\delta_k=\rho_k(0)\in\DD$, i.e., $\rho_k(z)=z\rho_{k-1}(z)+\delta_k\rho_{k-1}^*(z)$ for  $k\ge 1$.
All our spaces $\QQ_{n,k}$ will be subspaces of $\PP_n$, and when we write $\QQ_{n,k}^\perp$, we mean its orthogonal
complement in $\PP_n$.
We denote by $\varphi_a(z) = \frac{z +a}{1+\overline{a}z}$, the M\"obius Transform
associated with $a\in\CC\backslash\TT$.
We also assume that $\sum_{k=a}^b x_k$ is 0 and $\prod_{k=a}^b x_k$ is 1 if $b<a$.

\section{Quasi-paraorthogonal polynomials}\label{secQPOPUC}

As we have already said in the Introduction, quasi-orthogonal polynomials with respect to a measure $\mu$ on the real line have real coefficients, and as a consequence, their zeros are real or appear in complex conjugate pairs (symmetric with respect to $\RR$). The analogue on the unit circle, is that zeros are on $\TT$ or appear in symmetric pairs with respect to $\TT$.
That is: if $\alpha$ is a zero then also $1/\overline{\alpha}$ is a zero.
Note that this implies that if $P\in\PP_n\backslash\PP_{n-1}$, then $\alpha=0$ cannot be a zero of $P$. Thus, if the zeros are $\alpha_k\ne0$, $k=1,\ldots,n$, then we are looking for polynomials satisfying
\begin{eqnarray*}
P(z)&=&\nu\prod_{k=1}^n(z-\alpha_k)=
\nu\prod_{k=1}^n(z-\frac{1}{\overline{\alpha_k}})=
\nu c' z^n\prod_{k=1}^n\left(\overline{\frac{1}{\overline{z}}-{\alpha_k}}\right),\quad{c'}=
\prod_{k=1}^n(-{\alpha_k})\in\TT\\
&=&\frac{\nu}{\overline{\nu}}c' z^n\overline{P(\frac{1}{\overline{z}})}=\tau z^n\overline{P}(\frac{1}{z})=\tau P^*(z),\quad\tau=c'\frac{\nu}{\overline{\nu}}\in\TT.
\end{eqnarray*}
The notation $P^*(z)=z^n\overline{P(1/\overline{z})}$ is somewhat standard in the OPUC literature to denote the reciprocal of a polynomial $P$.
We have indeed that if $P(z)=\sum_{k=0}^n p_k z^k$ then $P^*(z)=\sum_{k=0}^n \overline{p_{n-k}}z^k$.
This explains the following definition.
\begin{definition}[{\bf invariant}]
We say that the polynomial $P\in\PP_n\backslash\PP_{n-1}$ is invariant
if and only if it satisfies $P=\tau P^*$, $\tau\in\TT$. The parameter $\tau$ is called the invariance parameter and $P$ is said to be $\tau$-invariant.
\end{definition}
\begin{remark}\label{rem2.4}
Note that if $P$ is a monic invariant polynomial, then its invariance parameter equals $P(0)\in\TT$ and if $\nu\in\CC\backslash\{0\}$, and $P$ is $\tau$-invariant, then $\nu P$ is $\tilde\tau$-invariant
with $\tilde\tau=\tau \nu/\overline{\nu}$.
\end{remark}

Thus, since we are interested in the symmetry of the zeros, we should only consider {\it invariant} polynomials.
Szeg\H{o} polynomials $\rho_n$ have all their zeros in $\DD$ and thus they are not invariant.
Quasi versions are obtained by imposing certain orthogonality conditions to subspaces of dimension $n-d$.
What should this subspace like for an invariant polynomial $Q_n\in\PP_n\backslash\PP_{n-1}$?

Let $Q_n$ be an invariant polynomial on $\TT$. Then
\[
\langle Q_n, z^k \rangle = \tau \langle Q^*_n, z^k \rangle = \tau \langle z^n \overline{Q_n}, z^k \rangle  = \tau \langle z^{n-k}, Q_n \rangle = \tau \overline{\langle Q_n, z^{n-k} \rangle},~~\tau\in\TT.
\]
Thus $\langle Q_n,z^{n-k}\rangle=0 \Leftrightarrow \langle Q_n,z^k\rangle=0$, for $0\le k\le n$.
An invariant polynomial on $\TT$ should always be orthogonal to a subspace that is spanned by powers of $z$ that are centrosymmetric in the sequence $(0,1,\ldots,n)$, that is, if it includes $z^k$, it should include $z^{n-k}$ as well.
Note that for QOPRL the number of orthogonality conditions is reduced by 1 as the order increases by 1.
However the symmetry in paraorthogonality implies that if we remove orthogonality to $z^{n-\ell}$ then we also remove orthogonality to $z^\ell$.
The number of orthogonality conditions, and thus also the order of QPOPUC, changes in steps of two.
So we define a nested sequence of subspaces of dimension $n-2\ell-1$ in $\PP_n$:
\[
\QQ_{n,2\ell+1}=\Span\{z^k: k=\ell+1,\ldots,n-\ell-1\},\quad 0\le\ell\le \left\lfloor \frac{n}{2}\right\rfloor-1.
\]
Invariant polynomials, orthogonal to this subspace for $\ell=0$
were coined {\it paraorthogonal} as they were introduced in \cite{JNT}, (see also \cite[Theorem~III]{Ger}), and they were later studied by many
\cite{DG91a,DG91b,CMV02,CMV06,CMNR16,CCBPP16,G,Wong,MISCFSS19,NS}.
We shall refer to them as quasi-paraorthogonal polynomials (QPOPUC) of order 1. More generally, we define the set of $(n,2\ell+1)$-QPOPUC of order $2\ell+1$ as all invariant polynomials of degree $n$ in
$\QQ^\perp_{n,2\ell+1}$.
The prefix quasi refers to the fact that they satisfy less paraorthogonality (i.e., symmetric orthogonality) conditions.

We recover  a representation of $(n,2\ell+1)$-QPOPUC in terms of the orthogonal polynomials obtained by Peherstorfer in \cite{Peh}.

\begin{theorem}\label{ThQuasi}
The monic $(n,2\ell+1)$-QPOPUC are given by
\[
Q_{n,2\ell+1}= zp_\ell \rho_{n-\ell-1}+\tau p_\ell^*\rho_{n-\ell-1}^*,~~p_\ell\in\PP_\ell\backslash \PP_{\ell-1} \,\, \mbox{monic and $Q_{n,2\ell+1}(0)=\tau\in\TT$.}
\]
The set of all monic $(n,2\ell+1)$-QPOPUC $Q_{n,2\ell+1}$ depends on $2\ell+1$ real free parameters.
\end{theorem}
\begin{proof}
We first note that
$\QQ_{n,2\ell+1}^\perp$ has dimension $n+1-(n-2\ell-1)=2\ell+2$.
Furthermore, $0=\langle z\rho_{n-\ell-1},z^k\rangle = \langle z^{t+1}\rho_{n-\ell-1},z^{k+t}\rangle=\langle \rho^*_{n-\ell-1},z^{k}\rangle = \langle z^t\rho^*_{n-\ell-1},z^{k+t}\rangle$, for $k=1,\ldots,n-\ell-1$, $t=0,\ldots,\ell$.
Thus  the $2\ell+2$ polynomials $\{z^{t+1}\rho_{n-\ell-1},z^t\rho_{n-\ell-1}^*\}_{t=0}^\ell$ are in $\QQ_{n,2\ell+1}^\perp$.
Moreover they are  independent because
$zp_\ell(z) \rho_{n-\ell-1}(z)+ q_\ell(z) \rho_{n-\ell-1}^*(z)\equiv0$ is only possible if
$p_\ell=q_\ell\equiv0$.
Indeed, if $p_\ell \not \equiv 0$, then  $\frac{z\rho_{n-\ell-1}}{\rho_{n-\ell-1}^*}\equiv-\frac{q_\ell}{p_\ell}$,
which is impossible because the rational function
$\frac{z\rho_{n-\ell-1}}{\rho_{n-\ell-1}^*}$ is irreducible of degree $n-\ell$ since $z\rho_{n-\ell-1}$ and $\rho_{n-\ell-1}^*$ have no common zeros
and the rational function $-\frac{q_\ell}{p_\ell}$ has degree at most $\ell<n-\ell$.
We get a similar contradiction if  $q_\ell \not \equiv 0$
so that $\QQ_{n,2\ell+1}^\perp=\{zp_\ell\rho_{n-\ell-1}+ q_\ell\rho_{n-\ell-1}^*: p_\ell,q_\ell\in\PP_\ell\}$.

Let $Q_n=z p_\ell\rho_{n-\ell-1}+q_\ell\rho_{n-\ell-1}^*$ be in $\QQ_{n,2\ell+1}^\perp$ with $p_\ell \in \PP_{\ell}\backslash\PP_{\ell-1}$ monic, then $Q_n$ is monic of degree $n$.
Because the representation w.r.t.\ a basis is unique, $Q_n$ is $\tau$-invariant if and only if $q_\ell= \tau p_\ell^*$. Therefore $Q_{n,2\ell+1}= z p_\ell \rho_{n-\ell-1}+\tau  p^*_\ell\rho_{n-\ell-1}^*$.
The $2\ell+1$ parameters represented by $\tau=Q_{n,2\ell+1}(0)\in\TT$ that depends on a real parameter and by the $\ell$ complex coefficients $p_\ell$, that correspond to $2\ell$ real parameters.
\end{proof}

The monic polynomial $p_\ell$ can be written as
\[
p_\ell(z)=\prod_{k=1}^\ell (z-\eta_k),\quad \eta_k\in\CC,~~k=1,\ldots,\ell.
\]
So we shall in the rest of the paper often use the $\ell$ zeros $\{\eta_k\}_{k=1}^\ell$ as parameters instead of the $\ell$ free coefficients of the monic polynomial $p_\ell$.

If $\sigma_{n,\ell}=\overline{p_\ell(0)}-\overline\tau{\delta_{n-\ell}}\ne0$, then $Q_n=Q_{n,2\ell+1}\not\in\QQ^\perp_{n,2\ell-1}$ (see Lemma~\ref{lemA} of the Appendix)
but there will be some `extra' orthogonality.
Because
$\langle Q_n,z^{\ell}\rangle=\tau\langle Q_n^*, z^\ell\rangle=\tau\overline{\langle Q_n,z^{n-\ell}\rangle}=\tau\overline{a_{n-\ell}}$
with $a_{n-\ell}=\langle Q_n,z^{n-\ell}\rangle\ne0$.
Define $\tilde\tau=\tau\frac{\sigma_{n,\ell}}{\overline{\sigma_{n,\ell}}}=\tau\frac{\overline{a_{n-\ell}}}{a_{n-\ell}}\in\TT$, then clearly $\langle Q_n,z^{n-\ell}-\tilde\tau z^\ell\rangle=0$.
Thus every $(n,2\ell+1)$-QPOPUC with $\sigma_{n,\ell}\ne0$ is orthogonal to
$\QQ_{n,2\ell+1}\oplus \Span\{z^{n-\ell}-\tilde\tau z^\ell\}$.
We shall refer to $\tilde\tau$ as the {\it orthogonality parameter} of $Q_{n,2\ell+1}$.

If $\sigma_{n,\ell}\ne0$ then the orthogonality parameter can be made explicit in an alternative expression
for $Q_{n,2\ell+1}$ as explained in Lemma~\ref{lemomega} of the Appendix.
If $\ell=0$, then $p_\ell\equiv1$ and thus $\sigma_{n,0}$ is always nonzero.

\section{Zeros of quasi-paraorthogonal polynomials}\label{secZero}

In the previous section we have described all monic and invariant quasi-paraorthogonal polynomials. We are now in a position to prove the main result that is an analogue on the unit circle of a very well known important result when dealing with measures supported on the real line that connects the number of orthogonality conditions with the number of zeros of odd multiplicity in the interior of the support of the measure.

It has been shown in the literature (see \cite{Sho,Jou,Mon,BCBVB10,BBMCQ13}) that an $(n,d)$-QOPRL associated with a measure whose support is an interval $[a,b]$ has at least $n-d$ zeros of odd multiplicity in $(a,b)$.
For $d=1$ we can use the parameter to fix one zero in $(a,b)$, for $d=2$, there are two parameters that can be used to fix $a$ and $b$ as nodes and for $d=3$
one may fix $a$, $b$ and some $\alpha\in(a,b)$, while in all these cases
it is guaranteed that there are $n$ simple zeros in the support $[a,b]\subset\RR$ of the measure.
These QOPRL for $d=1$ or $d=2$ are mainly derived in the context of Gauss-Radau and Gauss-Lobatto quadrature formulas respectively.

In this section we shall obtain similar results for QPOPUC assuming that the measure is supported on an arc $S\subset\TT$.

\subsection{A general statement}

The next Lemma gives a property of invariant polynomials that will be of interest for our purposes.
\begin{lemma}\label{invariante}
Let $P \in \PP_n \backslash \PP_{n-1}$ be an invariant polynomial. Then,
\begin{enumerate}
\item {\em Factoring structure on $\TT$}: if we denote by $\omega_1,\dots, \omega_m$ the zeros of $P$ of odd multiplicity on $\TT$, then
\begin{equation}\label{factoritation}
P(z) = K \cdot z^{r} \prod_{i=1}^m (z-\omega_i)\prod_{k=1}^{r}|z-\eta_k|^2, \quad K\neq 0, \quad z \in \TT,
\end{equation}
where $m+2r=n$ and $\eta_k$ are zeros of $P$ that may possibly coincide with $\omega_i$.
\item If $n$ is odd, then the number of zeros of odd multiplicity on $\TT$ is odd.
    If $n$ is even, then the number of zeros of odd multiplicity on $\TT$ is even (possibly zero).
\end{enumerate}
\end{lemma}
\begin{proof}
The second property is a direct consequence of the first one.
If $\alpha\not\in\TT$ is a zero of $P$, then $\frac{1}{\overline{\alpha}}= \frac{\alpha}{|\alpha|^2}$ is also a zero of $P$. So, the factorization of $P$ will have the factors
\begin{equation}\label{factor1}
(z-\alpha)\left(z-\frac{\alpha}{|\alpha|^2}\right)=-\frac{\alpha}{|\alpha|^2} z \left|z-\alpha\right|^2,\,\,z\in\mathbb{T}.
\end{equation}
If $\alpha$ is a zero of $P$ on $\TT$ and it has multiplicity at least two, then $P$ will have among its factors
\begin{equation}\label{factor2}
(z-\alpha)^2=-\alpha z|z-\alpha|^2,\,\,z\in\mathbb{T}.
\end{equation}
Thus, if $m$ is the number of zeros of odd multiplicity on $\TT$, then $P$ will have $m$ simple factors and the other factors are of the form (\ref{factor1})-(\ref{factor2}), yielding (\ref{factoritation}).
\end{proof}

We shall denote an oriented arc as an interval $[a,b]\subset\TT$ running counterclockwise from $a$ to $b$.
If $[a,b]\subset\TT$ denotes a (closed) arc, then the complementary arc with the same orientation is $(b,a)$. The square or rounded brackets are used as in the case of a real interval to indicate that boundary points are included or not.
Moreover, an arc can be indicated by three points like arc$(a,c,b)$, then this defines its orientation:
from $a$ over $c$ to $b$, and we write $a<c<b$.

\begin{theorem}\label{zeros}
For $n \geq 2$, $0 \leq \ell \leq \lfloor \frac{n}{2} \rfloor - 1$, and $\mathrm{supp}(\mu)=S=[a,b]\subset\TT$, an $(n,2\ell+1)$-QPOPUC $Q_{n,2\ell+1}$ has at least $n-2\ell$ zeros of odd multiplicity on the unit circle and $n-2\ell-1$ of them in $\disk{S}=(a,b)$.
\end{theorem}
\begin{proof}
Let us first consider the case $n=2s+1$.
Since a QPOPUC $Q_{n,2\ell+1}$ is invariant, we know from Lemma \ref{invariante} that $Q_{n,2\ell+1}$ has an odd number of zeros of odd multiplicity on $\TT$ and that it can be factored as
\[
Q_{n,2\ell+1}(z)=K z^{s-r}\prod_{j=1}^{t} (z-\omega_j)
\prod_{j=1}^u(z - \nu_j)
\prod_{j=1}^{s-r}|z-\eta_j|^2, \quad \quad K\neq 0, \quad \quad z=e^{i\theta} \in \TT
\]
with $t+u=2r+1$ odd.
We have distinguished between the $\eta_j$ with even multiplicity and for those of odd multiplicity we have the zeros $\nu_j$ that are in $\disk{S}$ and the distinct $\omega_i$ that are
not in $\disk{S}$, that is $\omega_i\in\TT\backslash \disk{S}=[b,a]$.
Without loss of generality we order the $\omega_i$ such that $b\le\omega_1<\omega_2<\cdots<\omega_t\le a$.
We define the arcs $S_i=[\omega_{2i-1},\omega_{2i}]$ with midpoints $\chi_i$ and $\psi_i(z) = \frac{z - \omega_{2i-1}}{z - \omega_{2i}}$
for $i=1,2,\ldots,\lfloor\frac{t}{2}\rfloor$.

Our aim is to prove that $u\ge n-2\ell-1$ which we do by contraposition, thus
suppose that $u \le n - 2\ell -2$.
We first treat the case where $t$ is even. Thus suppose $t=2k$, $k\in\NN$, then
\begin{eqnarray*}
\prod_{j=1}^k(z-\omega_{2j-1})\prod_{j=1}^k(z-\omega_{2j}) & = &  z^k \prod_{j=1}^k (-\omega_{2j}) \left[\prod_{j=1}^k(z-\omega_{2j-1})\prod_{j=1}^k(\overline{z}-\overline{\omega_{2j}})\right]\\[0.3cm]
 & = & z^k \prod_{j=1}^k (-\omega_{2j}) \left[ \prod_{j=1}^k\psi_j(z) \prod_{j=1}^k |z-\omega_{2j}|^2\right],\quad z \in \TT.
\end{eqnarray*}
Furthermore, if
\[
\Theta(z) = K \left(\prod_{j=1}^k \frac{-\omega_{2j}}{\overline{\psi_j(\chi_j)}}\right) z^{ s - r +k} \prod_{j=1}^u (z - \nu_j),
\]
then
\[
Q_{n,2\ell+1}(z)\overline{\Theta(z)}=|K|^2\prod_{j=1}^{s-r}|z-\eta_j|^2\prod_{j=1}^u|z-\nu_j|^2\prod_{j=1}^k|z-\omega_{2j}|^2
\prod_{j=1}^k\frac{\psi_j(z)}{\psi_j(\chi_j)}, \quad z \in \TT.
\]
By Lemma \ref{signo},
\[
z  \in \disk{S} \Rightarrow \frac{\psi_j(z)}{\psi_j(\chi_j)} < 0 ,\quad j=1, \dots, k,
\]
so that $Q_{n,2\ell+1}(z)\overline{\Theta(z)}$, has a constant sign in $\disk{S}$.

But $\langle Q_{n,2\ell+1}, \Theta \rangle = 0$ because $\Theta \in \QQ_{n,2\ell+1}$. Indeed, since $\Theta \in \Span\left\{ z^{s-r+k},\ldots,z^{s-r+k+u} \right\}$, it only remains to prove that $s-r+k \geq \ell+1$ and that $s-r+k+u \leq n-\ell-1$. The assumptions are $n=2s+1$, $t+u=2r+1$, $u \leq n-2\ell-2$ and $t=2k$ (notice that $u$ is odd). So,
\[
s-r+k = s-r + \frac{2r+1-u}{2} = s + \frac{1}{2} - \frac{u}{2} \geq \ell + 1
\]
 and
\[
 s - r + k + u = s + \frac{1}{2} + \frac{u}{2} \leq s + \frac{1}{2} + \frac{n-2\ell-2}{2} = n - \ell-1.
\]
This proves that $Q_{n,2\ell+1}$ has at least $n-2\ell-1=2(s-\ell)$ zeros of odd multiplicity in $\disk{S}=(a,b)$, and since $n$ is odd, $Q_{n,2\ell+1}$ has $n-2\ell$ zeros on the unit circle.

Now, we consider $t= 2k+1$, then $u$ is even and $u \leq 2(s-\ell) -2$.
We consider
\[
\Theta(z) = K \left(\prod_{j=1}^k \frac{-\omega_{2j}}{\overline{\psi_j(\chi_j)}}\right) z^{ s - r +k} (z - \omega_{2k+1})\prod_{j=1}^u (z - \nu_j),
\]
then
\[
Q_{n,2\ell+1}(z)\overline{\Theta(z)}=|K|^2\prod_{j=1}^{s-r}|z-\eta_j|^2\prod_{j=1}^u|z-\nu_j|^2\prod_{j=1}^k|z-w_{2j}|^2 |z-\omega_{2k+1}|^2
\prod_{j=1}^k\frac{\psi_j(z)}{\psi_j(\chi_j)},\quad z \in \TT.
\]

By the same argument, $Q_{n,2\ell+1}(z)\overline{\Theta(z)}$ has a constant sign in $\disk{S}$, which is a contradiction since $\Theta \in \QQ_{n,2\ell+1}$ and thus $\langle Q_{n,2\ell+1}, \Theta \rangle = 0$.
Indeed, $\Theta \in \Span\left\{ z^{s-r+k},\ldots,z^{s-r+k+u+1} \right\}$ and now, $2k+1+u = 2r+1$ ($u$ is even) and we have supposed that $u \leq 2(s - \ell)- 1$, hence  $u \leq 2(s - \ell- 1)$.
So, $s-r+ k = s- \frac{u}{2} \geq \ell + 1$ and $s - r + k + u + 1 = s + 1 + \frac{u}{2} \leq s + 1 + (s-\ell -1) = n - \ell-1$.
Again in this case $u \geq 2(s-\ell)$. This proves that $Q_{n,2\ell+1}$ has at least $n-2\ell-1=2(s-\ell)$ zeros of odd multiplicity in $\disk{S}=(a,b)$, and since $n$ is odd, $Q_{n,2\ell+1}$ has $n-2\ell$ zeros on the unit circle.

When $n$ is even the proof is similar.
\end{proof}

\subsection{$Q_{n,1}$ and one prescribed zero}
If we can employ the invariance parameter $\tau$ (i.e.\ $Q_{n,1}(0)$ for a monic $Q_{n,1}$) to place an extra zero of odd multiplicity on $\disk{S}$, then
this would imply that  the $n$ simple zeros of $Q_{n,1}$ are in $\disk {S}$.
We shall study the situation using Blaschke products.
Since the zeros $z_{k,n}$ of the Szeg\H{o} polynomial $\rho_n$ lie in $\DD$,
we can define the Blaschke products of degree $n$ (with zeros in $\DD$) as follows
\begin{equation}\label{FHdef}
F_n(z)=z\frac{\rho_{n-1}(z)}{\rho_{n-1}^{*}(z)}=\prod_{s=0}^{n-1} \varphi_{z_{s,n}}(z)
~~\text{where}~~ z_{0,n} = 0 ~~ \mbox{and} ~~ \rho_{k}(z)=\prod_{l=1}^{k}\left( z-z_{l,k} \right), ~~ k \in \NN.
\end{equation}
By the Argument Principle, $F_n(z)$ goes around the unit circle exactly $n$ times when $z$ wraps around the origin (a $n$-to-1 map), and since all their zeros are in $\DD$ the map $z\mapsto F_n(z)$ preserves the orientation. If in addition the support of the measure is an arc $S$, as it is in our case, the following result is a consequence of the previous theorem.
\begin{corollary}\label{Blaschke}
Let $\mu$ be a positive measure supported on an arc $S\subset\TT$ with associated orthogonal polynomials $\{\rho_k\}_{k\in\NN}$ and let $F_n$ be the Blaschke product as defined above in (\ref{FHdef}).
Then for given $t\in \TT$, $F_n(z)=t$  will have at least $n-1$ solutions in $\disk{S}$.
\end{corollary}

\begin{theorem}\label{zero}
Let $\mu$ be a positive measure supported on $S=[a,b]\subsetneq\TT$, and $n\ge 2$.
Consider the QPOPUC $Q_{n,1}(z)=z\rho_{n-1}(z)+\tau\rho_{n-1}^*(z)$ and define
$\tau_a=F_n(a)$ and $\tau_b=F_n(b)$ with $F_n$ as in (\ref{FHdef}).
Then
$Q_{n,1}$ has all its simple zeros in $\disk{S}$ (or in $S$) if and only if $-\tau \in (\tau_a,\tau_b)$ (or $-\tau\in [\tau_a,\tau_b]$ respectively). Moreover $\alpha\in\TT$ is one of these zeros  if and only if $-\tau=F_n(\alpha)$.
\end{theorem}

\begin{proof}
The zeros $\{z_i\}_{i=1}^n$ of $Q_{n,1}$ are the roots of the equation $F_n (z)= -\tau$.
We denote the (ordered) roots of $F_n(z) =  \tau_a$ by $\{\alpha_i\}_{i=1}^n$ and the (ordered) roots of $F_n(z) = \tau_b$, by $\{\beta_i\}_{i=1}^n$, with $\alpha_1=a$ and $\beta_n=b$. By the properties of $F_n$, $z_i \in (\alpha_i,\beta_i)$, $i =1,2, \ldots, n$.
\end{proof}

\begin{remark}\label{rem01}
Note also that the only prefixed zero $\alpha\in\disk{S}$ that is allowed must be in one of the $n$ intervals $(\alpha_i,\beta_i)$, $i=1,\ldots, n$.
Compare with \cite[Theorem~2.9]{BCBVB10} for QOPRL.
\end{remark}

\subsection{$Q_{n,3}$ and 2 or 3 prescribed zeros}

A monic QPOPUC
\begin{equation}\label{eqQn3}
Q_{n,3}(z)=z (z-\eta)\rho_{n-2}(z)+ \tau (1-\overline{\eta}z)\rho_{n-2}^*(z)
\end{equation}
of order 3, has at least $n-2$ zeros of odd multiplicity on $\TT$, $n-3$ of them are in interior of the arc where the measure is supported, and it depends on three real parameters that can be used to fix up to three zeros, ensuring that all the zeros of $Q_{n,3}$ are simple and placed on the support of measure.
Let us start by fixing 2 nodes.

\begin{figure}[!ht]
\begin{center}
\includegraphics[scale=0.6]{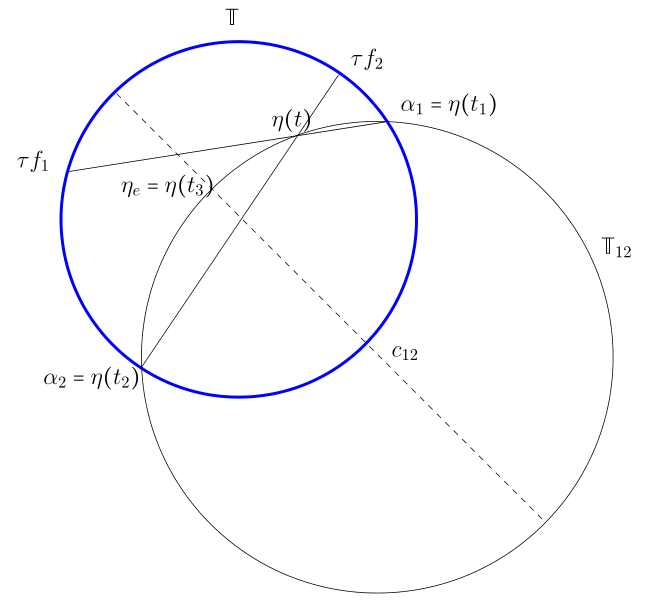}
\caption{\label{fig0}The unit circle $\mathbb{T}$ and the circle $\mathbb{T}_{12}$ with $\mathbb{T}\cap\mathbb{T}_{12}=\{\alpha_1,\alpha_2\}$ and the two secant lines through $(\alpha_1,{\tau}f_1)$ and $(\alpha_2,{\tau}f_2)$ intersecting at $\eta(\tau)$.}
\end{center}
\end{figure}

If $\eta \in \DD$, then $B_n(z) = \frac{z (z-\eta)\rho_{n-2}(z)}{(1-\overline \eta z) \rho_{n-2}^*(z)}$ is a Blaschke product of degree $n$ and $B_n(z) = \tau \in \TT$ has $n$ simple solutions on $\TT$.
The conditions $Q_{n,3}(\alpha_1) = Q_{n,3}(\alpha_2) = 0$ lead to
\begin{equation}\label{eqfi0}
\alpha_i(\alpha_i-\eta) \rho_{n-2}(\alpha_i) + \tau (1-\overline\eta \alpha_i) \rho_{n-2}^*(\alpha_i) = 0, \quad i =1,2,
\end{equation}
or equivalently,
\begin{equation}\label{eqfi}
(\alpha_i - \eta) + {\tau} {f_i} (1-\overline \eta \alpha_i) = 0, \quad f_i = \overline{F_{n-1}(\alpha_i)}, \quad  i =1,2.
\end{equation}
This gives the system
\begin{equation}\label{eqfi2}
\left[\begin{array}{cc}
1 & \tau {f_1}\alpha_1\\
1 & \tau {f_2}\alpha_2
\end{array}\right]
\left[\begin{array}{c}\eta\\\overline{\eta}\end{array}\right]=
\left[\begin{array}{c}\tau{f_1}+\alpha_1\\
\tau{f_2}+\alpha_2\end{array}\right].
\end{equation}
If $f_1\alpha_1\ne f_2\alpha_2$, this represents two secants for the circle $\TT$
\begin{equation}\label{secants}
 \begin{array}{c}
\eta + {\tau} f_1 \alpha_1 \overline{\eta} = {\tau} f_1 + \alpha_1
\quad\text{and}\quad
\eta + {\tau} f_2 \alpha_2 \overline{\eta}= {\tau} f_2 + \alpha_2,~~\eta\in\CC.
\end{array}
\end{equation}
The first secant passes through $\alpha_1$ and ${\tau} f_1$ and the second  passes through $\alpha_2$ and ${\tau} f_2$.
If $f_1 = f_2$, then $\eta = {\tau} f_1\in \TT$, and this should be excluded because $\eta$ must be in $\DD$.
These secants intersect in exactly one finite point if and only if $\alpha_1 f_1 \neq \alpha_2 f_2$ and  the intersection point is
\begin{equation}\label{eta}
\eta=c_{12}+ {\tau}a_{12},
~~\overline{c_{12}}=\frac{f_1-f_2}{f_1\alpha_1-f_2\alpha_2},~~\overline{a_{12}}=\frac{\alpha_1-\alpha_2}{f_1\alpha_1-f_2\alpha_2}.
\end{equation}
By Lemma~\ref{chords} of the Appendix, $\eta \in \DD$ if and only if the arcs $(\alpha_1, {\tau}f_1, \alpha_2)$ and $(\alpha_1, {\tau}f_1, {\tau}f_2)$ have opposite orientation. (See also Figure~\ref{fig0}.)
Note that we can only define $(\alpha_1, {\tau}f_1, \alpha_2)$ and $(\alpha_1, {\tau}f_1, {\tau}f_2)$ as arcs if $\tau f_1\ne\alpha_2$ and $\tau f_2\ne\alpha_1$ which is implicitly assumed also in the next Theorem.

\begin{theorem}\label{pts2}
Let $n\ge3$.
Given $\tau\in\TT$, and distinct $\{\alpha_i\}_{i=1}^2\subset \TT$.
Define $\{f_i\}_{i=1}^2\subset\TT$ as in (\ref{eqfi}).
Assume the $f_i$ are distinct and $\alpha_1 f_1\ne \alpha_2 f_2$.
Then the monic $(n,3)$-QPOPUC $Q_{n,3}$ of (\ref{eqQn3}) is uniquely defined and has all its zeros simple and on $\TT$, $\alpha_1$ and $\alpha_2$ being two of them if $\tau$ is such that
the arcs $(\alpha_1, {\tau}f_1, \alpha_2)$ and $(\alpha_1, {\tau}f_1, {\tau}f_2)$ have opposite orientation. The
parameter $\eta$ in (\ref{eqQn3}) belongs to $\DD$ and depends on $\tau$ as given by (\ref{eta}).
\end{theorem}
\begin{remark}
Note that there are infinitely many $\tau$ that allow us to fix $\alpha_1$ and $\alpha_2$ ensuring that $Q_{n,3}$  has all its zeros simple on $\TT$.
\end{remark}

A further geometric interpretation can be obtained as follows (see also \cite{JaRei}).
Because $\tau\in\TT$, the identity (\ref{eta})  says that $\eta\in \TT_{12}$ which is the circle with center $c_{12}$
and with radius $r_{12}=|a_{12}|$, i.e.\ (see also Figure~\ref{fig0})
\[ 
\TT_{12}=\{ \eta(t)=c_{12}+t\cdot r_{12}: t\in\TT\}.
\] 
The $\eta$ from (\ref{eta}) corresponds to $t={\tau}{t_0}$ where $t_0={a_{12}}/|a_{12}|$.
For $t_1=\alpha_1\overline{f_2}t_0$ we get $\eta(t_1)=\alpha_1$ and for $t_2=\alpha_2\overline{f_1}t_0$
we get $\eta(t_2)=\alpha_2$.
Thus $\TT\cap\TT_{12}=\{\alpha_1,\alpha_2\}$.
Since there are two different points of intersection, this means that there are values for $t\in\TT$ that will deliver $\eta(t)\in\DD$
and others for which $\eta(t)\in\EE$.
It is clear that the line connecting the origin with $c_{12}$, i.e., the mediatrix of the chord $[\alpha_1,\alpha_2]$, will intersect $\TT$ at the points $\pm\sqrt{\alpha_1\alpha_2}$ and the circle $\TT_{12}$ at two points:
one in $\DD$ that is closest to the origin and one in $\EE$ that is the farthest from the origin.
The former is given by $\eta_e=c_{12}(1-\frac{r_{12}}{|c_{12}|})$, which is $\eta(t_3)$ with $t_3=-\frac{c_{12}}{|c_{12}|}$.
Thus $\eta(t)\in\TT_{12}\cap\DD$
if and only if $t$ is in the interior of the $\arc(t_1,t_3,t_2)\subset\TT$
for which $\eta(t)$ describes the $\arc(\eta(t_1),\eta(t_3),\eta(t_2))=\arc(\alpha_1,\eta_e,\alpha_2)=\TT_{12}\cap\DD$.
Taking into account the value of $t_0$, $t\in\arc(t_1,t_3,t_2)$ if and only if
$\tau\in\arc\left(\alpha_1\overline{f_2},-\frac{\overline{f_1}-\overline{f_2}}{|f_1-f_2|}
\frac{\alpha_1-\alpha_2}{|\alpha_1-\alpha_2|},\alpha_2\overline{f_1}\right)$.

All this assumes $f_1\alpha_1\ne f_2\alpha_2$, i.e., the two secants have a
single finite intersection point.
If $f_1\alpha_1= f_2\alpha_2$, then the secants are parallel or coincide.
\begin{enumerate}
\item $(\alpha_1-\alpha_2)+\tau(f_1-f_2)\ne0$, then the two equations of (\ref{secants}) are contradictory. The two secants are parallel and intersect at $\eta=\infty$.
\item $(\alpha_1-\alpha_2)+\tau(f_1-f_2)=0$, then the two equations of the system (\ref{eqfi2}) are identical and reduce to $\eta+\overline{\eta}\alpha_1\alpha_2=\alpha_1+\alpha_2$, the secant through $\alpha_1=\tau f_2$ and $\alpha_2=\tau f_1$ (See Figure~\ref{fig0}).
In this case there are infinitely many solutions for $\eta$ given by $\eta=\alpha_1+t(\alpha_2-\alpha_1)$ with $\eta\in\DD$ if $0<t<1$.
This line is in fact a degenerate case of the circle $\TT_{12}$ with an infinite radius and a center at infinity.
The mediatrix of the chord $[\alpha_1,\alpha_2]$ intersects $\TT_{12}$ at the center $\eta_e=(\alpha_1+\alpha_2)/2$ and at $\infty$. Since $\eta$ has to be in $\DD$, it should be on the degenerate $\arc(\alpha_1,\eta_e,\alpha_2)=\TT_{12}\cap\DD$ which is the open interval $(\alpha_1,\alpha_2)$.
\end{enumerate}

\begin{remark}\label{remtau1}
Note that $Q_{n,3}$ has 3 real parameter incorporated by $\eta\in\CC$ and $\tau\in\TT$. In the generic case ($\alpha_1 f_1\ne\alpha_2 f_2$) we can, depending on a chosen $\tau$ fix the zeros $\alpha_1$ and $\alpha_2$ by selecting an appropriate $\eta$. However, to guarantee that all $n$ zeros are simple and on $\TT$, we had to restrict $\tau$ to a subarc of $\TT$.
In the degenerate case (2) above, $\alpha_1$ and $\alpha_2$ are chosen such that $\alpha_1\overline{f_2}=\alpha_2\overline{f_1}=\tau$, which means that now $\tau$ is fixed, but it leaves
$\eta$ still free to vary on the secant line through $\alpha_1$ and $\alpha_2$, which corresponds to one real parameter $t$.
To guarantee also in this case that all $n$ zeros are simple and on $\TT$,
we should restrict $\eta$ to the subinterval $(\alpha_1,\alpha_2)$, i.e.\ $0<t<1$. The latter can be seen as an arc of the degenerate circle $\TT_{12}$ through $\alpha_1$ and $\alpha_2$ and center at $\infty$.
\end{remark}

With this analysis we can turn the existence Theorem~\ref{pts2} into an explicit characterization of $\tau$ as follows.
\begin{theorem}\label{pts2b}
Let $n\ge3$.
Given $\tau\in\TT$, and $\{(\alpha_i,f_i)\}_{i=1}^2\subset \TT\times \TT$, satisfying the conditions of Theorem~\ref{pts2}.
Then the monic $(n,3)$-QPOPUC $Q_{n,3}$ as in (\ref{eqQn3})
with $n$ simple zeros on $\TT$, $\alpha_1$ and $\alpha_2$ being two of them,
is uniquely defined if
$\tau\in\arc\left(\alpha_1\overline{f_2},
-\frac{\overline{f_1}-\overline{f_2}}{|f_1-f_2|}\frac{\alpha_1-\alpha_2}{|\alpha_1-\alpha_2|},
\alpha_2\overline{f_1}\right)$.
The parameter $\eta$ of (\ref{eqQn3}) belongs to $\DD$ and depends on $\tau$ as given by (\ref{eta}).
\end{theorem}

The previous argument shows that if the two prefixed zeros allow a solution that guarantees $n$ simple zeros on $\TT$ then we should choose free parameter
($\tau$ or $t$) such that $\eta$ is on some subarc in $\DD$.
This remaining freedom will allow us to fix a third zero under certain conditions as we shall currently show.

The conditions $Q_{n,3}(\alpha_1)=Q_{n,3}(\alpha_2)=Q_{n,3}(\alpha_3)=0$ lead to
\[
\alpha_i(\alpha_i-{\eta})\rho_{n-2}(\alpha_i)+ \tau(1-\overline\eta \alpha_i)\rho_{n-2}^*(\alpha_i)=0, \quad i =1,2,3
\]
or equivalently,
\[
\overline{\tau} (\alpha_i-\eta) + f_i (1 - \alpha_i \overline \eta) = 0, \quad f_i=\overline{ F_{n-1}(\alpha_i)},~~ \quad i = 1,2,3,
\]
where $\eta \in \CC$ and $\tau \in \TT$.
If $1-\alpha_i \overline{\eta} \neq 0$, $i=1,2$, we can reformulate the system as
\begin{equation}\label{autom}
\overline{\tau}\varphi_{-\eta}(\alpha_i)=-f_i, \quad i =1,2,3.
\end{equation}
Then, $Q_{n,3}$ has $n$ simple zeros if the solution of the above interpolation problem is an automorphism of the unit disk.
It well know that there exists a unique automorphism of the unit disk if and only if $\arc(\alpha_1,\alpha_2,\alpha_3)$ and  $\arc(f_1,f_2,f_3)$ have the same orientation, \cite[p.~397]{Saff}, see also \cite[p.~280 and p.~296, exercise~30]{R}, and that it is a particular case of a more general situation: see \cite[Theorem~4.1 and Example~5.1]{B}.

\begin{theorem}\label{3pointsfixed}
Let $n\ge3$.
Given distinct $\{\alpha_i\}_{i=1}^3 \subset \TT$,
and assume $f_i=\overline{F_{n-1}(\alpha_i)}$, $i=1,2,3$ with $F_{n-1}$ as in (\ref{FHdef}) are also distinct.
Then there exists a unique monic $(n,3)$-QPOPUC $Q_{n,3}$ as in (\ref{eqQn3}), having all its zeros simple on $\TT$, with $\alpha_1$, $\alpha_2$ and $\alpha_3$ being three of them  if $\arc(\alpha_1,\alpha_2,\alpha_3)$ and $\arc(f_1,f_2,f_3)$ have the same orientation.
The required $\eta$ of (\ref{eqQn3}) belongs to $\DD$.
\end{theorem}

If the support of $\mu$ is an arc $S=[a,b]\subsetneq\TT$, then  Theorem~\ref{zeros} guarantees that $Q_{n,3}$  has $n-2$ zeros of odd multiplicity on $\TT$ of which only $n-3$ are in the interior of the support.
We can however prove that {\it all} the zeros are in the support if we choose $\alpha_1=a$ and $\alpha_2=b$ in the following way.

The interpolation conditions (\ref{eqfi}) become for $\eta \not\in\{a,b\}$,
\begin{equation}\label{eqint}
\varphi_{-\eta}(a)=-\tau f_a~~\text{and}~~
\varphi_{-\eta}(b)=-\tau f_b~~\text{where $f_a=\overline{F_{n-1}(a)}$ and $f_b=\overline{F_{n-1}(b)}$.}
\end{equation}
We know this has a unique solution $\eta$ for each $\tau\in\TT$ if and only if $f_a \neq f_b$ and $a f_a \neq b f_b$
and the latter will always be true in this case as a consequence of Corollary \ref{Blaschke}.
For this $\eta$, $Q_{n,3}$ has zeros at $a$ and $b$, and thus $Q_{n,3}(z)=(z-a)(z-b)P_{n-2}(z)$ with $P_{n-2}\in\PP_{n-2}\backslash \PP_{n-3}$ monic.
These boundary zeros $a$ and $b$ are in addition to the $(n-3)$ zeros that $Q_{n,3}$ has in $(a,b)$ according to Theorem~\ref{zeros}.
As $Q_{n,3}$ is $\tau$-invariant and orthogonal to $\QQ_{n,3}$,
$P_{n-2}$ is $\hat\tau$-invariant with $\hat\tau=\overline{ab}\tau$ and orthogonal to $\QQ_{n-2,1}$ with respect to
\begin{equation}\label{eqmuhat}
d\hat\mu(z)=\sqrt{\overline{ab}}(z-a)(z - b)\overline{z} d\mu(z),~~ z \in S.
\end{equation}
We assume the square root is chosen such that this measure $\hat\mu$ is positive on $\disk{S}$, which is possible because with $z=e^{i\theta}$, $\sqrt{a}=e^{i\theta_a/2}$, and $\sqrt{b}=e^{i\theta_b/2}$,
\[
\sqrt{\overline{ab}}\; e^{-i \theta}(e^{i\theta}- e^{i \theta_a})(e^{i\theta}- e^{i \theta_b}) = 4\sin\Bigl(\frac{\theta- \theta_a}{2}\Bigr)\sin\Bigl(\frac{\theta_b - \theta}{2}\Bigr)> 0, \quad \forall \theta \in (\theta_a,\theta_b).
\]
Thus, $P_{n-2}$ is a monic $\hat\tau$-invariant, $(n-2,1)$-QPOPUC for the positive measure $\hat\mu$.
Let $\{\hat{\rho}_n\}_{n\in \NN}$ be the sequence of monic orthogonal polynomials associated with $\hat\mu$, then $P_{n-2}=z \hat\rho_{n-3}+\hat\tau \hat\rho_{n-3}^*$.

Applying Theorem~\ref{zero} thus gives us the following result.

\begin{theorem}\label{thm23}
Let $\mu$ be a measure supported on an arc $S=[a,b]\subsetneq\TT$
and define the associated positive measure $\hat\mu$ as in (\ref{eqmuhat}) with monic orthogonal polynomial sequence $\{\hat{\rho}_n\}_{n\in \NN}$.
For $k\in\NN$, let $F_k$ be defined as in (\ref{FHdef}) and  $\hat{F}_k(z)={z
\hat\rho_{k-1}(z)}/{\hat\rho_{k-1}^*(z)}$.
Assume $n\ge3$ and define
$f_a = \overline{F_{n-1}(a)}$, $f_b=\overline{F_{n-1}(b)}$
$\hat{\tau}_a=\hat{F}_{n-2}(a)$, $\hat{\tau}_b=\hat{F}_{n-2}(b)$, and $\hat\tau=\overline{ab}\tau$.
Furthermore set $\eta = c_{12}+ a b \hat\tau a_{12}$, where
$\overline{c_{12}}=\frac{f_a-f_b}{af_a-bf_b}$, $\overline{a_{12}}=\frac{a-b}{af_a -bf_b}$ and $\hat\tau\in\TT$.
Finally, let $Q_{n,3}$ as in (\ref{eqQn3})
be a monic $(n,3)$-QPOPUC.
Then
\begin{enumerate}
\item
$Q_{n,3}$ has all its simple zeros on $S$ where  $a$ and $b$ are two of them if and only if
$-\hat\tau\in(\hat\tau_a,\hat\tau_b)$
\item
$Q_{n,3}$ has all its simple zeros in $S$ where  $a,b$ and $\alpha\in(a,b)$  are three of them if and only if $-\hat\tau=\hat{F}_{n-2}(\alpha)\in(\hat\tau_a,\hat\tau_b)$.
\end{enumerate}
\end{theorem}

\subsection{$Q_{n,2\ell+1}$ with $2\ell$ or $2\ell+1$ prescribed zeros}

For the more general case where $\ell>1$, the analysis becomes very hard
when the measure is supported on a subarc of $\TT$. So for the rest of this section we shall assume that the measure is supported on the whole circle $\TT$.

Let us first fix $2\ell$ zeros for the $(n,2\ell+1)$-QPOPUC
\begin{equation}\label{Qn2lp1}
Q_{n,2\ell+1}(z)=zP_\ell(z)\rho_{n-\ell-1}(z)+\tau P_\ell^*(z)\rho_{n-\ell-1}^*(z),\quad P_\ell(z)=\prod_{k=1}^\ell (z-\eta_k),~~\eta_k\in\CC.
\end{equation}
From Theorem~\ref{zeros} we know it has at least $(n-2\ell)$ zeros of odd multiplicity on $\TT$.
It remains to verify which conditions on $\tau\in \TT$ and $\{\eta_i\}_{i=1}^\ell$ must be satisfied to allow us to fix $2\ell$ nodes, while ensuring that all the zeros of $Q_{n,2\ell+1}$ are simple and on $\TT$.

If $\eta_i \in \DD$, $i=1,\ldots,\ell$, then $B_n(z) = \frac{z P_\ell(z)\rho_{n-\ell-1}(z)}{P_\ell^*(z) \rho_{n-\ell-1}^*(z)}$ is a Blaschke product of degree $n$ with all zeros in $\DD$ and thus $B_n(z) = \tau \in \TT$ has $n$ simple solutions on $\TT$.
We already know that the conditions $Q_{n,2\ell+1}(\alpha_j)=0$, $j=1,\ldots,2\ell$, lead to
\[
\alpha_i P_\ell(\alpha_i)
 \rho_{n-\ell-1}(\alpha_i) + \tau P_\ell^*(\alpha_i)
 \rho_{n-\ell-1}^*(\alpha_i)= 0, \quad i = 1, \ldots, 2\ell,
\]
or equivalently, with
$f_{i}=\overline{F_{n-\ell}(\alpha_i)}$, $i = 1,\ldots,2\ell$, ($F_{n-\ell}$ as in (\ref{FHdef})):
\begin{equation}\label{efsl2}
  P_\ell(\alpha_i) +  \tau f_i P_\ell^*(\alpha_i)= 0, \quad  i = 1,\ldots,2\ell.
\end{equation}
Obviously, a possible solution will depend on the selection of the nodes $\{\alpha_i\}$ and how the corresponding $\{f_i\}$ behave.
This problem has a matrix interpretation, that extends the case $\ell=1$.

Let us assume that there exists a unique solution $P_\ell$ for the system (\ref{efsl2}) and that it is given by
$P_\ell(z)= p_0 + p_1 z + \cdots + p_{\ell-1} z^{\ell-1} + z^\ell$.
Set
$\mathcal{V}_1=[\alpha_i^j]_{i=1,\ldots,\ell}^{j=0,\ldots,\ell-1}$ and
$\mathcal{V}_2=[\alpha_i^j]_{i=\ell+1,\ldots,2\ell}^{j=0,\ldots,\ell-1}$
the Vandermonde matrices associated with
$\{ \alpha_i\}_{i=1}^{\ell}$ and $\{ \alpha_j\}_{j=\ell+1}^{2\ell}$, respectively, and furthermore set
\[
 \begin{array}{ll}
  \v{f}_1 = [{f_1},{f_2},\dots,{f_\ell}]^T,& \mathcal{F}_1=\mathrm{diag}\left(\v{f}_1 \right),\\
  \v{f}_2 = [{f_{\ell+1}},{f_{\ell+2}},\dots, {f_{2\ell}}]^T,& \mathcal{F}_2=\mathrm{diag}\left(\v{f}_2 \right),\\
  \v{d}_1 = [\alpha^\ell_1,\alpha^\ell_2,\dots, \alpha^\ell_\ell]^T,&  \mathcal{D}_1=\mathrm{diag}\left( \v{d}_1 \right),\\
  \v{d}_2 = [\alpha^\ell_{\ell+1},\alpha^\ell_{\ell+2},\dots, \alpha^\ell_{2\ell}]^T,& \mathcal{D}_2=\mathrm{diag}\left( \v{d}_2 \right),
 \end{array}
\]
then, from (\ref{efsl2})
\[
 p_0+p_1\alpha_i+\cdots+p_{\ell-1}\alpha_i^{\ell-1} +\tau{f_i}\alpha_i^\ell\left(\overline{p}_{\ell-1} \overline{\alpha_i^{\ell-1}} +\cdots + \overline{p}_0 \right)=-\tau{f_i}-\alpha_i^\ell, \quad i=1,\ldots,2\ell.
\]
Setting $\v{p}=\left[ p_0, \;p_1, \;\cdots, \;p_{\ell-1} \right]^T$, these equations are equivalent to
\[
\mathcal{V}_k \v{p} + \tau{\mathcal{F}_k} \mathcal{D}_k \overline{\mathcal{V}_k}\overline{\v{p}} = -\tau{\v{f}_k} - \v{d}_k, \quad k=1,2,
\]
or in matrix form, the vector $[\v{p}^T,\overline{\v{p}}^T]^T$ is the unique solution of
\begin{equation}\label{system2l}
\left[ \begin{array}{cc}
\mathcal{V}_1 & \tau{\mathcal{F}_1}\mathcal{D}_1\overline{\mathcal{V}_1} \\
\mathcal{V}_2 & \tau{\mathcal{F}_2}\mathcal{D}_2\overline{\mathcal{V}_2}
\end{array}\right]
\left[ \begin{array}{c}\v{p}\\ \overline{\v{p}}\end{array}\right]=
\left[ \begin{array}{c}-\tau\v{f}_1-\v{d}_1 \\ -\tau{\v{f}_2} - \v{d}_2 \end{array}\right].
\end{equation}
The uniqueness of the solution implies that the coefficient matrix is nonsingular.

\begin{theorem}\label{lemaadd}
Let $n\ge2\ell+1$. Given $\tau\in\TT$ and distinct $\{\alpha_i\}_{i=1}^{2\ell}\subset\TT$,
and assume $f_i=\overline{F_{n-\ell}(\alpha_i)}$, $i=1,\ldots,2\ell$ with $F_{n-\ell}$ defined in (\ref{FHdef}) are distinct and satisfy
\begin{equation}\label{conditionaddendum1}
\prod_{1 \leq i < j \leq \ell}
\frac{
\left( \alpha_{\ell+i}-\alpha_{\ell+j} \right)
\left( \overline{\alpha_i} - \overline{\alpha_j} \right)
}{
\left( \overline{\alpha_{\ell+i}}-\overline{\alpha_{\ell+j}}\right)
\left(\alpha_i-\alpha_j\right)
}
\neq \prod_{i=1}^{\ell} \frac{{f_{\ell+i}}\alpha_{\ell+i}^\ell}{{f_i}\alpha_i^\ell}.
\end{equation}
Then the monic $\tau$-invariant $(n,2\ell+1)$-QPOPUC $Q_{n,2\ell+1}(z)=zP_\ell(z)\rho_{n-\ell-1}(z)+\tau P_\ell^*(z)\rho_{n-\ell-1}^*(z)$ with $n$ simple zeros on $\TT$, and $\alpha_1, \alpha_2,\dots, \alpha_{2\ell}$ being $2\ell$ of them is uniquely defined if
all $s_k(0) \in \DD$, $k=1,\dots, \ell$, where
 \begin{equation}\label{Schur}
s_\ell(z)= \frac{P_\ell(z)}{P_\ell^*(z)}~~\text{and}~~s_{k-1}(z)=\frac{1}{z}\frac{s_k(z)-s_k(0)}{1-\overline{s_k(0)}s_k(z)},~~k=\ell,\ell-1,\ldots,2,
\end{equation}
and $P_\ell(z)=p_0+ p_1z + \cdots + p_{\ell-1}z^{\ell-1}+z^\ell$ where
$\v{p}=\left[ p_0 \;p_1 \;\cdots \;p_{\ell-1} \right]^T$  is the unique solution of the system (\ref{system2l}).
The vector $\v{p}$ is given by
\begin{equation}\label{conditionaddendumP}
 \v{p}=-\left[
 \left(\overline{\mathcal{V}_1}\right)^{-1}
 \overline{\mathcal{D}_1} \overline{\mathcal{F}_1} \mathcal{V}_1 -
 \left(\overline{\mathcal{V}_2}\right)^{-1}
 \overline{\mathcal{D}_2} \overline{\mathcal{F}_2} \mathcal{V}_2 \right]^{-1} \cdot
 \left[
 \left(\overline{\mathcal{V}_1}\right)^{-1}\left( \tau\overline{\v{d}_1}+\overline{\v{f}_1} \right) -
 \left(\overline{\mathcal{V}_2}\right)^{-1}\left( \tau\overline{\v{d}_2}+\overline{\v{f}_2} \right)
 \right],
\end{equation}
where we use the notation that was introduced above.
\end{theorem}
\begin{proof}
We consider the following auxiliary problem. Let  $P_\ell(z)=p_0+ p_1z + \cdots + p_{\ell-1}z^{\ell-1}+z^\ell$ and
$Q_\ell(z)=1+q_1z+ \cdots+q_{\ell-1}z^{\ell-1}+q_\ell z^\ell$ verify
\[
 p_0+p_1\alpha_i+\cdots+p_{\ell-1}\alpha_i^{\ell-1} +\tau{f_i}\alpha_i^\ell\left(q_\ell+q_{\ell-1}\overline{\alpha_i}+\cdots +q_1\overline{\alpha_i}^{\ell-1} \right)=-\tau{f_i}-\alpha_i^\ell, \quad i=1,\ldots,2\ell.
\]
Setting $\v{p}=\left[ p_0, \;p_1, \;\cdots, \;p_{\ell-1} \right]^T$ and
$\v{q}=\left[ q_1, \;q_2, \;\cdots, \;q_\ell \right]^T$, these are equivalent to
\[
\mathcal{V}_k \v{p} + \tau{\mathcal{F}_k} \mathcal{D}_k \overline{\mathcal{V}_k}\mathcal{S}_{\ell}\v{q} = -\tau {\v{f}_k} - \v{d}_k, \quad k=1,2,
\]
where $\mathcal{S}_\ell=[\delta_{i,\ell+1-j}]_{i,j=1}^\ell$ is the $\ell\times\ell$ flipping matrix which maps a vector upside down.
Thus
\begin{equation}\label{conditionaddendum2}
\left[ \begin{array}{cc}
\mathcal{V}_1 & \tau{\mathcal{F}_1}\mathcal{D}_1\overline{\mathcal{V}_1} \\
\mathcal{V}_2 & \tau{\mathcal{F}_2}\mathcal{D}_2\overline{\mathcal{V}_2}
\end{array}\right]
\left[ \begin{array}{c}\v{p}\\ \mathcal{S}_{\ell}\v{q} \end{array}\right]=
\left[ \begin{array}{c}-\tau\v{f}_1-\v{d}_1 \\ -\tau{\v{f}_2} - \v{d}_2 \end{array}\right].
\end{equation}
The condition (\ref{conditionaddendum1}) ensures that the associated matrix is not singular, hence the system has a unique solution.

That $P_\ell=Q_\ell^*$ can be seen as follows.
Multiply by
$\overline{\tau}\mathrm{diag}(\overline{\mathcal{F}_1} \overline{\mathcal{D}_1},
               \overline{\mathcal{F}_2} \overline{\mathcal{D}_2})$
from the left in (\ref{conditionaddendum2}) to obtain
\[
\left[ \begin{array}{cc}
\overline{\tau\mathcal{F}_1}\overline{\mathcal{D}_1}\mathcal{V}_1 & \overline{\mathcal{V}_1} \\
\overline{\tau\mathcal{F}_2}\overline{\mathcal{D}_2}\mathcal{V}_2 & \overline{\mathcal{V}_2} \end{array}\right]
\left[ \begin{array}{c} \v{p} \\ \mathcal{S}_{\ell}\v{q} \end{array}\right]=
\left[ \begin{array}{c} -\overline{\tau\v{d}_1}-\overline{\v{f}_1} \\ -\overline{\tau\v{d}_2}-\overline{\v{f}_2} \end{array}\right].
\]
By taking the complex conjugate and interchanging the block columns we get
\[
\left[ \begin{array}{cc}
\mathcal{V}_1 & \tau{\mathcal{F}_1}\mathcal{D}_1\overline{\mathcal{V}_1} \\
\mathcal{V}_2 & \tau{\mathcal{F}_2}\mathcal{D}_2\overline{\mathcal{V}_2}  \end{array}\right]
\left[ \begin{array}{c} \mathcal{S}_{\ell} \overline{\v{q}} \\ \overline{\v{p}} \end{array}\right]=
\left[ \begin{array}{c} -{\tau\v{f}_1}-\v{d}_1 \\ -{\tau\v{f}_2} - \v{d}_2 \end{array}\right].
\]
Hence, from the unicity of the solution, $\v{p}=\mathcal{S}_{\ell}\overline{\v{q}}$, and so, $P_\ell=Q_{\ell}^*$.

For the second part, observe that (\ref{conditionaddendumP}) follows directly from (\ref{conditionaddendum2}) by block Gaussian elimination.

The third part described in (\ref{Schur}) is the classical Schur-Cohn test (see \cite{Schur,Cohn}) for the polynomial $P_\ell$.
\end{proof}
Note that the yeast of this proof is essentially the same as in the Theorem~\ref{pts2}.
We require that the system expressing the interpolation conditions has a unique solution, and require that the solution $P_\ell$ has its zeros in $\DD$, which was obtained in Theorem~\ref{pts2} and by restricting $\tau$ to an arc, which is here replaced by an a Schur-Cohn test.

\begin{remark}\label{remtaul}
We point out that (\ref{conditionaddendum1}) does not depend on $\tau$, so neither does the existence of $Q_{n,2\ell+1}$.
Whether or not $Q_{n,2\ell+1}$ in Theorem~\ref{lemaadd} will have the desired properties
depends on the value of $\tau\in\TT$.
If (\ref{conditionaddendum1}) is satisfied, there will exist a polynomial $P_\ell$ such that the interpolation  conditions are satisfied.
However, whether or not the other conditions of the previous theorem can be satisfied, in particular the conditions that the polynomial $P_{\ell}$ will have zeros in $\DD$, depends continuously on the parameter $\tau\in\TT$.
All the $s_k(0)$ will be in $\DD$ for $\tau$ on certain arcs of $\TT$.
If $\tau$ leaves such an arc, this means that some Schur parameter $s_{\ell_0}(0)$, and thus also a zero $\eta$ of $P_\ell$, leaves $\DD$ by crossing $\TT$. At that moment $\eta$ becomes one of the zeros of $Q_{n,2\ell+1}$.
For $\tau$ outside these arcs, $Q_{n,2\ell+1}$ can have zeros of multiplicity larger than 1 and/or pairs of zeros
$(z_i,1/\overline{z_i})$ with $z_i\in\DD$.
To find out what exactly are the pathological cases requires further analysis of the location and multiplicity of the zeros, which can be obtained by a deeper exploration of the behaviour of the Schur parameters. See for example
\cite{Bis,VM} and similar discussions in the literature.
However, keeping in mind Theorem~\ref{zeros}, we can always assure that there are at least $n-2\ell$ zeros on $\TT$ with odd multiplicity.
\end{remark}

The freedom to choose $\tau$ in the indicated arcs of $\TT$ may allow us to fix another additional node. More precisely, if we want
to fix $2\ell+1$ zeros $\{\alpha_i\}_{j=1}^{2\ell+1}\subset\TT$ of $Q_{n,2\ell+1}$, it should be verified that
\[
\alpha_i P_\ell(\alpha_i) \rho_{n-\ell-1}(\alpha_i) +\tau P_\ell^*(\alpha_i)\rho_{n-\ell-1}^*(\alpha_i)=0,\quad j=1,2,\dots,2\ell+1.
\]
Equivalently, if $f_i=\overline{F_{n-\ell}(\alpha_i)}\in \TT$ and $\lambda= \sqrt{\overline{\tau}}$,
\begin{equation}\label{efsl2p1}
\lambda P_\ell(\alpha_i)+f_i \overline\lambda P^*_\ell(\alpha_i)=0,\quad i=1,2,\dots,2\ell+1.
\end{equation}
Also this  problem has a matrix interpretation where the main difference with the previous case is that now $\tau=Q_{n,2\ell+1}(0)=\overline{\lambda^2}$ is a parameter that we have to determine by solving the system.

More precisely, if we denote
$P_\ell(z)=\prod_{j=1}^\ell(z - \eta_j) = p_0 + p_1 z + \cdots + p_{\ell-1} z^{\ell-1} + z^\ell$,
then from (\ref{efsl2p1})
\[
\lambda \left( p_0 + p_1\alpha_i +\cdots + p_{\ell-1}\alpha_i^{\ell-1}+\alpha_i^\ell\right) +{f_i}\overline{\lambda} \left(1 +\overline{p}_{\ell-1} \alpha_i +\cdots + \overline{p}_0 \alpha_i^{\ell} \right)=0, \quad i=1,\ldots,2\ell+1.
\]
Setting $\v{p}=\lambda[ p_0,p_1,\ldots,p_{\ell-1},1 ]^T$,
and $\mathcal{V}=[\alpha_i^j]_{i=1,\ldots,2\ell+1}^{j=0,\ldots,\ell}$ ,  $\v{f} = [{f_1},{f_2},\dots,{f_{2\ell+1}}]^T$ and $ \mathcal{F}=\mathrm{diag}\left(\v{f} \right)$,
this is equivalent to the homogeneous system
\[
\mathcal{V} \v{p} + \mathcal{F} \mathcal{V} \v{p^*} = 0,~~\text{or}~~
\left[ \mathcal{V}~~\mathcal{FV}\right] \left[\begin{array}{c}\v{p}\\\v{p}^*\end{array}\right]=\v{0},
\]
and the vector $[\v{p}^T,\v{p}^{*T}]^T$ is the unique solution, up to a multiplicative constant means that the matrix $[\mathcal{V}~~\mathcal{FV}]\in\CC^{ (2\ell+1)\times(2\ell+2)}$ has full rank $2\ell+1$. Since $P_\ell$ has degree $\ell$ deleting column $\ell+1$ in this matrix makes it square and nonsingular.

\begin{theorem}\label{lemaadd2}
Let $n\ge2\ell+1$. Given distinct numbers $\{\alpha_i\}_{i=1}^{2\ell+1} \subset \TT$, and let $f_i=\overline{F_{n-\ell}(\alpha_i)}$, $i=1,\dots, 2\ell+1$
with $F_{n-\ell}$ as defined in (\ref{FHdef}) be also distinct.
Denote by $\mathcal{V}_1=[\alpha_i^j]_{i=1,\ldots,\ell}^{j=0,\ldots,\ell-1}$ and
$\mathcal{V}_2=[\alpha_i^j]_{i=\ell+1,\ldots,2\ell+1}^{j=0,\ldots,\ell}$
the square Vandermonde matrices associated with $\{ \alpha_i\}_{i=1}^{\ell}$ and $\{ \alpha_i\}_{i=\ell+1}^{2\ell+1}$, respectively, and the rectangular Vandermonde matrices
$\mathcal{V}'_1=[\alpha_i^j]_{i=1,\ldots,\ell}^{j=0,\ldots,\ell}\in\CC^{\ell\times(\ell+1)}$ and
$\mathcal{V}'_2=[\alpha_i^j]_{i=\ell+1,\ldots,2\ell+1}^{j=0,\ldots,\ell-1}\inf\CC^{(\ell+1)\times \ell}$,
and furthermore set
\[
 \begin{array}{ll}
  \v{f}_1 = [{f_1},{f_2},\dots,{f_\ell}]^T,& \mathcal{F}_1=\mathrm{diag}\left(\v{f}_1 \right),\\
  \v{f}_2 = [{f_{\ell+1}},{f_{\ell+2}},\dots, {f_{2\ell+1}}]^T,& \mathcal{F}_2=\mathrm{diag}\left(\v{f}_2 \right),\\
  \v{d}_1 = [\alpha^\ell_1,\alpha^\ell_2,\dots, \alpha^\ell_\ell]^T, &
  \v{d}_2 = [\alpha^\ell_{\ell+1},\alpha^\ell_{\ell+2},\dots, \alpha^\ell_{2\ell+1}]^T,
  \end{array}
\]
and assume
\begin{equation}\label{regular2l+1}
\Delta =\det \left[ \begin{array}{cc}
\mathcal{V}_1 & \mathcal{F}_1\mathcal{V}'_1 \\
\mathcal{V}'_2 & \mathcal{F}_2\mathcal{V}_2
\end{array}\right]\ne0.
\end{equation}
Then there exists a unique monic invariant $(n,2\ell+1)$-quasi-paraorthogonal polynomial $Q_{n,2\ell+1}(z)=zP_\ell(z)\rho_{n-\ell-1}(z)+\tilde\tau P_\ell^*(z)\rho_{n-\ell-1}^*(z)$
with  all its zeros simple on $\TT$, with $\alpha_1$, $\alpha_2,\dots$, $\alpha_{2\ell+1}$  being $2\ell+1$  of them
if
all $s_k(0) \in \DD$, $k=1,\dots, \ell$, where
\begin{equation}\label{Schur2}
s_\ell(z)= \frac{P_\ell(z)}{P_\ell^*(z)}~~\text{and}~~
s_{k-1}(z)=\frac{1}{z}\frac{s_k(z)-s_k(0)}{1-\overline{s_k(0)}s_k(z)},~~k=\ell,\ell-1,\ldots,2,
\end{equation}
where $P_\ell(z) = \sum_{k=1}^{\ell-1} p_k z^k+z^\ell$ is given by $\v{p}=[ p_0,p_1,\ldots,p_{\ell-1} ]^T$ such that
\begin{equation}\label{solution2lp1}
 \v{p}=\left[
 \mathcal{V}_1 - \mathcal{F}_1\mathcal{V}'_1 \mathcal{V}_2^{-1}\overline{\mathcal{F}_2} {\mathcal{V}'_2} \right]^{-1} \cdot \left[\v{d}_1- \mathcal{F}_1\mathcal{V}'_1 \mathcal{V}_2^{-1}\overline{\mathcal{F}_2}\v{d}_2 \right].
 \end{equation}
\end{theorem}
\begin{proof}
The proof is analogous to the previous case. We consider the following auxiliary problem. Let  $\tilde{P}_\ell(z)=\tilde{p}_0+ \tilde{p}_1z + \cdots + \tilde{p}_{\ell-1}z^{\ell-1}+\tilde{p}_\ell z^\ell$ with $\tilde{p}_\ell\ne0$ and
$Q_\ell(z)=q_0+q_1z+ \cdots+q_{\ell-1}z^{\ell-1}+q_\ell z^\ell$ satisfy
\[
 \tilde{p}_0+\tilde{p}_1\alpha_i+\cdots+ \tilde{p}_\ell \alpha_i^\ell +{f_i}\left(q_0 + q_1 \alpha_i + \cdots + q_{\ell}\alpha_i ^{\ell}\right)= 0,\quad i=1,\ldots,2\ell+1.
\]
It is easy to verify that $Q_\ell = \tilde{P}_\ell^*$,
then setting $\tilde{\v{p}}=\left[ \tilde{p}_0,\tilde{p}_1,\ldots,\tilde{p}_{\ell-1} \right]^T$ and
$\v{q}=\left[ \overline{\tilde{p}_\ell}, \overline{\tilde{p}_{\ell-1}},\ldots, \overline{\tilde{p}_0} \right]^T$, the above equations are equivalent to the linear system
\[
\left[ \begin{array}{cc}
\mathcal{V}_1 &  \mathcal{F}_1\mathcal{V}'_1 \\
\mathcal{V}'_2 & \mathcal{F}_2\mathcal{V}_2
\end{array}\right]\left[\begin{array}{c}\tilde{\v{p}} \\ \v{q} \end{array}\right] = -\tilde{p}_\ell \left[\begin{array}{c} \v{d}_1 \\ \v{d}_2 \end{array}\right],\;\; \tilde{p}_\ell \neq 0
\]
which has a unique solution by (\ref{regular2l+1}). Then, (\ref{solution2lp1}) follows directly from  block Gaussian elimination taking into account that $\mathcal{V}_1$ and $ \mathcal{F}_2\mathcal{V}_2$ are nonsingular and
\[
\frac{\overline{\tilde{p}_\ell}}{\tilde{p}_\ell} = -\det \left[\begin{array}{cc}
\mathcal{V}'_1 & \mathcal{D}_1\mathcal{F}_1\mathcal{V}_1 \\
\mathcal{V}_2 &  \mathcal{D}_2\mathcal{F}_2\mathcal{V}'_2\end{array}\right]/
\det\left[\begin{array}{cc}
\mathcal{V}_1 & \mathcal{F}_1\mathcal{V}'_1 \\
\mathcal{V}'_2 & \mathcal{F}_2\mathcal{V}_2
\end{array}\right] = -\frac{1}{\Delta}\det \left[\begin{array}{cc}
\mathcal{V}'_1 & \mathcal{D}_1\mathcal{F}_1\mathcal{V}_1 \\
\mathcal{V}_2 &  \mathcal{D}_2\mathcal{F}_2\mathcal{V}'_2\end{array}\right],
\]
where 
$\mathcal{D}_1=\mathrm{diag}( \alpha_1,\ldots,\alpha_\ell)$,  $\mathcal{D}_2=\mathrm{diag}( \alpha_{\ell+1},\ldots,\alpha_{2\ell+1})$.
Since all the entries in the matrices are on $\TT$, we can rewrite the determinant of the numerator as
\[
\det \left[\begin{array}{cc}
\mathcal{D}_1^\ell\mathcal{F}_1 & \\ & \mathcal{D}_2^\ell\mathcal{F}_2
\end{array}\right]
\det \left[\begin{array}{cc}
\overline{\mathcal{F}_1}\overline{\mathcal{V}_1} & \overline{\mathcal{V}'_1} \\
\overline{\mathcal{F}_2}\overline{\mathcal{V}'_2} &  \overline{\mathcal{V}_2}\end{array}\right]=(-1)^{\ell}\,\overline\Delta
\prod_{j=1}^{2\ell+1} f_j \alpha_j^\ell ~~
\text{thus}~~\frac{\overline{\,\tilde{p}_\ell}\,}{\tilde{p}_\ell} =(-1)^{\ell+1}  \frac{\,\overline\Delta\,}{\Delta} \prod_{j=1}^{2\ell+1} f_j \alpha_j^\ell.
\]
Moreover, if  $s_\ell(z)=P_\ell(z)/P_\ell^*(z)$, with $P_\ell(z)=
\tilde{P}_\ell(z)/\tilde{p}_\ell=\sum_{k=1}^{\ell-1} ({\tilde{p}_k}/{\tilde{p}_\ell}) z^k + z^{\ell}$, then  $P_\ell$ will have all its zeros in $\DD$ if and only if all the $s_k(0)\in\DD$ for $k=\ell,\ldots,1$ by the Schur-Cohn test.
Setting $P_\ell(z)=\prod_{i=1}^\ell (z-\eta_i)$ with $\eta_i\in\DD$, $i=1,\ldots,\ell$, it is clear that
\[
Q_{n,2\ell+1}(z)=z\prod_{i=1}^\ell (z-\eta_i)\rho_{n-\ell-1}(z)+\frac{\;\overline{\tilde{p}_\ell}\;}{\tilde{p}_\ell}
\prod_{i=1}^\ell (1-\overline{\eta_i}z)\rho_{n-\ell-1}^*(z)
\]
is the monic invariant $(n,2\ell+1)$-QPOPUC with invariance parameter $\tau=Q_{n,2\ell+1}(0)=\frac{\;\overline{\tilde{p}_\ell}\;}{\tilde{p}_\ell}$ and  all its zeros simple on $\TT$, with $\alpha_1, \alpha_2,\dots, \alpha_{2\ell+1}$  being $2\ell+1$  of them.
\end{proof}
\begin{remark}\label{remtau2}
Note that prefixing $2\ell+1$ zeros will in general define $Q_{n,2\ell+1}$ uniquely, which includes its invariance parameter $\tau$.
Thus there is no `varying' $\tau$.
Generically, it is impossible to fix more than $2\ell+1$ prefixed zeros, unless we select a zero $\alpha_{2\ell+2}$ that happens to be one of the zeros of $Q_{n,2\ell+1}$ that is already uniquely defined by the first $2\ell+1$ zeros.
A discussion somewhat similar to Remark~\ref{remtaul} could be made here if we choose $\{\alpha_i\}_{i=1}^{2\ell}$ fixed and consider $\alpha_{2\ell+1}$ as a varying parameter defining $\tau$, but we shall not elaborate on that.
\end{remark}

We note also that, with a different normalization for $P_\ell$, a similar technique was used in \cite[Theorem~2.2, Corollary~2.3]{gla2} to find Blaschke product interpolants of degree $\ell$ for $2\ell+1$ points on $\TT$.

\section{Applications to quadrature formulas on the unit circle}\label{secQF}

In this section we shall be concerned with the numerical estimation of integrals of the form
\begin{equation}\label{integral}
I(f)=\int f(z)d\mu(z),
\end{equation}
$f$ being a Riemann-Stieltjes integrable function with respect to $\mu$, by means of quadrature formulas (q.f.) of the form
\begin{equation}\label{qf}
I_n(f)=\sum_{s=1}^n \lambda_s f(z_s),
\end{equation}
with distinct nodes $\{ z_s \}_{s=1}^{n}\subset \supp(\mu)$ and
weights defined by $\lambda_s=I(L_{n,s})$ with $\{L_{n,s}\}_{s=1}^n$ the Lagrange Laurent polynomial interpolation basis for the nodes $\{z_s\}_{s=1}^n$ in an $n$ dimensional subspace of the form $\LL_{-p,q}$, $q+p=n-1$, with
$p,\ell\in\NN$ such that $0<2\ell+1<n$ and $\ell\le p\le n-\ell-1$.
This q.f.\ is obviously exact in $\LL_{-p,q}$.
If all weights in (\ref{qf}) are positive we call $I_n(f)$ a positive quadrature formula (positive q.f.).
Let us denote the nodal polynomial by $W_n(z)=\prod_{s=1}^n (z-z_s)$. Note that it is $W_n(0)$-invariant.

It is well known that the $(n,d)$-QOPRL with $d$ prefixed zeros can be used as the nodal polynomials of Gauss-type formulas, i.e., positive quadrature formulas with maximal domain of exactness like Gauss-Radau ($d=1$) or Gauss-Lobatto ($d=2$) or more general quadrature formulas with prescribed nodes (\cite{BCBVB10}).

The purpose of this section is to apply $(n,2\ell+1)$-QPOPUC as nodal polynomials to obtain similar results for Szeg\H{o}-type q.f.\
namely to characterize positive quadrature formulas with maximal domain of validity when some of the nodes are prefixed. This will correspond to q.f.\ whose nodal polynomial is a $(n,2\ell+1)$-QPOPUC with $2\ell$ or $2\ell+1$ prefixed nodes.

\subsection{Szeg\H{o}-type quadrature}

As a consequence of the density of Laurent polynomials in the space of continuous functions defined on an arc $S\subseteq\TT$ with respect to the uniform norm, q.f.\ like (\ref{qf}) are usually constructed by imposing
$I(L)=I_n(L)$ for all $L\in \Span\{ z^{-p},\ldots,z^q \}=\LL_{-p,q}$, $p,q\in \NN$ with $p+q$ as large as possible
(in this case, we say that $I_n(f)$ is exact in $\LL_{-p,q}$ and it has at least degree of exactness $p+q$).
The {\em Szeg\H{o} q.f.}\ introduced in \cite{JNT}  are positive q.f.\ whose degree of exactness is $2n-2$.
The nodal polynomial is $W_n=Q_{n,1}$, a
monic invariant $(n,1)$-QPOPUC $Q_{n,1}=z\rho_{n-1}+\tau\rho_{n-1}^*=\nu_{n,1}^{-1}[\rho_n+\tilde\tau\rho_n^*]$ in which the invariance parameter $\tau =Q_{n,1}(0)$ or, equivalently, the orthogonality parameter $\tilde\tau$, is free.
It is exact in $\LL_{-(n-1),n-1}$, that is a subspace of Laurent polynomials of dimension $2n-1$, despite the fact that the q.f.\ depends on $2n$ parameters ($n$ nodes and $n$ weights, all varying with $\tau$).
Recall that Gaussian q.f.\ also depend on $2n$ parameters but they are exact in $\PP_{2n-1}$, a space of dimension $2n$.
So we may expect that there must also be some space of dimension $2n$ in which  the Szeg\H{o} quadrature is exact which corresponds to the orthogonality represented by the orthogonality parameter $\tilde\tau$. This was first observed in \cite{NS} and \cite{Bul2}.
The authors in \cite{NS}, found that the Szeg\H{o} q.f.\  also integrated exactly
$z^n-\frac{\tilde\tau\tau}{z^{n}}$ thus achieving a subspace of dimension $2n$.
So we shall consider subspaces of Laurent polynomials that are invariant under the involution $f(z)\mapsto f_*(z)=\overline{f(1/\overline{z})}$ and we shall say that a q.f.\ on $\TT$ of the form (\ref{qf}) has degree of exactness $e$ if it is exact for such an invariant subspace of dimension $e$.
This and other issues in the construction of q.f.\ on the unit circle, were studied in a more general context by some of the authors in \cite{CDP}.
In that paper, it was proved that the resulting q.f.\ always have real weights $\{ \lambda_k \}_{k=1}^{n}$, at least $\lfloor \frac{e}{2}\rfloor + 1 $ of them are positive where $e\ge n-1$ is the degree of exactness of the q.f.\ \cite[Theorem~2.3]{CDP} and there always exist two positive q.f.\ that are exact in a maximal domain of dimension $2n$ \cite[Corollary~2.9]{CDP} when $\supp(\mu)=\TT$.
Also an analogue of a Jacobi-type theorem (see Theorem~\ref{Jacobi} below) was obtained \cite[Theorem~2.6]{CDP}.

This idea of constructing a q.f.\ whose nodal polynomial is a QPOPUC is briefly summarized as follows.
As in the real line situation, we need to define a nested sequence $\{{\Ll}_n\}_{n\in \NN}$ of subspaces of Laurent polynomials, such that ${\Ll}_{0}=\Span\{ 1 \}$, ${\Ll}_{n-1} \subset {\Ll}_n$, $\dim \left({\Ll}_{n}\right)=n+1$ and $\bigcup_{n\in \NN}{\Ll}_n={\LL}$.
This can be done by starting from a sequence $\{\omega_k\}_{k \geq 1} \subset \TT$ and defining
\begin{equation}\label{spaces}
{\Ll}_{2k}:=\LL_{-k,k}, \quad
\quad {\Ll}_{2k+1}={\Ll}_{2k+1}(\omega_{k+1}):={\Ll}_{2k}\oplus \Span\left\{z^{k+1}- \frac{\omega_{k+1}}{z^{k+1}}\right\}, \quad k \geq 0.
\end{equation}
Note that $\Ll_n$ is invariant under the involution $L\mapsto L_*$,
that is $L(z) \in {\Ll}_n \Leftrightarrow \overline{L(1/\overline{z})} \in {\Ll}_n$.
The next result  gives the paraorthogonality conditions that must be satisfied by the nodal polynomial of a q.f.\ $I_n(f)$ of the form (\ref{qf}) to be exact in subspaces ${\Ll}_{n-1+k}$, $0 \leq k \leq n$ of the form (\ref{spaces}) and not in $\Ll_{n+k}$.
It was proved in \cite{NS} for $k=n-1$ and in the context of Laurent polynomials in \cite[Theorem 2.6]{CDP} for any $k$.
For completeness we include a proof using the notation and the definitions introduced in Section~\ref{secQPOPUC} inspired by the results  obtained in \cite{NS} for $k=n-1$.
The result can be stated as follows:

\begin{theorem}\label{Jacobi}
Let $n,\ell,p \in \NN$ such that $0 < 2\ell+1 < n$, $\ell \leq p \leq n - \ell-1$. Let $I_n(f)$ be a q.f.\ of the form (\ref{qf}) exact in $\mathbb{L}_{-p,n-1-p}$. Let $W_n$ be its nodal polynomial and $a_{n-\ell}=\langle W_n,z^{n-\ell}\rangle\ne0$.
Then
\begin{enumerate}
\item $I_n(f)$ is exact in ${\Ll}_{2(n-\ell-1)}$ if and only if
$W_n$ is a monic $(n,2\ell+1)$-QPOPUC.
\item $I_n(f)$ is exact in ${\Ll}_{2(n-\ell)-1}(\omega)$ with $\omega= \tilde\tau W_n(0)$
if and only if $W_n$ is an $(n,2\ell+1)$-QPOPUC and $\tilde\tau=\tau\frac{\overline{a_{n-\ell}}}{a_{n-\ell}}$.
\end{enumerate}
\end{theorem}
\begin{proof}
Let $I_n(f)$ be exact in $\LL_{-p,q}$, and denote by $L_{n-1}(f,x)$ the Laurent polynomial interpolating $f$ in $\LL_{-p,n-1-p}$ in the nodes $\{x_i\}_{i=1}^n$ of the q.f. Let $L \in\Ll_{2(n-\ell-1)}$, then
\[
L(z) - L_{n-1}(L,z) = \frac{W_n(z) P(z)}{z^{n-\ell-1}}, \quad P \in \PP_{n-2\ell-2},
\]
integrating the equality,
\[
I(L) = I_n(L) , \, \forall \, L \, \in {\cal L}_{2(n-\ell-1)}  \quad \Leftrightarrow \quad 0 =  I(\frac{W_n P}{z^{n-\ell-1}}) = \langle W_n, Q \rangle,\, \forall  Q \in \mathbb{Q}_{n,2\ell+1}.
\]
It only remains to prove that additionally $z^{n-\ell}-\omega z^{-(n-\ell)}$ will be integrated exactly
if and only if $\omega=\tilde\tau\tau$.

Therefore we proceed as follows.
Assume that $Q_n:=Q_{n,2\ell+1}$ is a monic $(n,2\ell+1)$-QPOPUC with $\langle Q_n,z^{n-\ell}\rangle\ne0$ then its orthogonality parameter $\tilde\tau$ is well defined. If $\tau=Q_n(0)$, then
\begin{equation}\label{eqA1}
z^{-\ell}[z^n-Q_n]\in z^{-\ell}\LL_{0,n-1}=\LL_{-\ell,n-\ell-1}
~~\text{and}~~
z^{-(n-\ell)}[1-\overline{\tau}Q_n]\in z^{-(n-\ell)}\LL_{1,n}=\LL_{-(n-\ell-1),\ell}
\end{equation}
and because both $\LL_{-\ell,n-\ell-1}$ and $\LL_{-(n-\ell-1),\ell}$ are subspaces of $\LL_{-(n-\ell-1),n-\ell-1}$, (recall that $2\ell+1< n$) the quadrature will be exact for these Laurent polynomials:
\[
[I-I_n]\Bigl(z^{n-\ell}-z^{-\ell}Q_n\Bigr)=0~~\text{and}~~
[I-I_n]\Bigl(z^{-(n-\ell)}-\overline{\tau}z^{-(n-\ell)}Q_n\Bigr)=0
\]
but because $Q_n$ is the nodal polynomial, $I_n(z^{-\ell}Q_n)=0=I_n(z^{-(n-\ell)}Q_n)$, so that
\[
I(z^{n-\ell})=I_n(z^{n-\ell})+I(Q_n z^{-\ell}) ~~\text{and}~~
I(z^{-(n-\ell)})=I_n(z^{-(n-\ell)})+\overline{\tau}I(Q_n z^{-(n-\ell)}).
\]
Thus for any $\omega\in\TT$
\[
I\left(z^{n-\ell}-\frac{\omega}{z^{n-\ell}}\right)=I_n\left(z^{n-\ell}-\frac{\omega}{z^{n-\ell}}\right)+
\langle Q_n,z^{\ell}\rangle-\omega
\overline{\tau}\langle Q_n,z^{n-\ell}\rangle.
\]
We know that $I(Q_n z^{-(n-\ell)})=\langle Q_n,z^{n-\ell}\rangle\ne0$  and because $Q_n$ is $\tilde\tau$-orthogonal, $\langle Q_n,z^{\ell}\rangle=\tilde\tau \langle Q_n,z^{n-\ell}\rangle$ so that
\[
\langle Q_n,z^{\ell}\rangle-\omega\overline{\tau}\langle Q_n,z^{n-\ell}\rangle=
(\tilde\tau-\omega\overline{\tau})\langle Q_n,z^{n-\ell}\rangle=0
\Leftrightarrow
\omega=\tau\tilde\tau.
\]
\end{proof}

This theorem states a one-to-one correspondence between all QPOPUC that have  all their zeros simple and in the support of the measure and all q.f.\ exact in a certain subspace $\Ll_m$.

states that all QPOPUC that have  all their  single zeros on the support of the measure, are nodal polynomials of a q.f., exact in a certain subspace $\Ll_m$.
We have seen in previous sections that if $\supp(\mu)$ is an arc $S\subsetneq\TT$, then
$Q_{n,2\ell+1}$ has at least $n-2\ell$ zeros of odd multiplicity on $\TT$ and $n-2\ell-1$ of them in $\disk{S}$. It depends on $2\ell+1$ free real parameters that can be used to preselect up to $2\ell+1$ zeros simple and on $S$.
It remains to analyse when this will deliver positive q.f. with maximal domain of exactness.

\subsection{Szeg\H{o} and Szeg\H{o}-Radau quadrature}

Since $\nu_{n,1}=1+\tilde\tau\overline{\delta_n}\ne0$, it follows from Lemma~\ref{lemomega} in the Appendix, that $Q_{n,1}(z)=z\rho_{n-1}+\tau\rho_{n-1}^*(z)=\nu_{n,1}^{-1}[\rho_n(z)+\tilde\tau\rho_{n-1}^*(z)]$
and
\[
\tau=Q_{n,1}(0)=\tilde\tau\frac{\overline{\nu_{n,1}}}{\nu_{n,1}}=\varphi_{\delta_n}(\tilde\tau)
\Leftrightarrow \tilde\tau=\varphi_{-\delta_n}(\tau).
\]

We know that $Q_{n,1}$ has always $n$ simple zeros on $\TT$.
It was shown in \cite{JNT} that Szeg\H{o}'s q.f.\ has $Q_{n,1}$ as its nodal polynomial, it is exact in $\Ll_{2n-2}=\LL_{-(n-1),n-1}$ and its weights are positive. By \cite{NS} or Theorem~\ref{Jacobi}, we know that it is also exact in $\Ll_{2n-1}(\omega)$ with $\omega=\tilde\tau\tau=\tilde\tau\varphi_{\delta_n}(\tilde\tau)=\tau\varphi_{-\delta_n}(\tau)$.
The freedom to pick an arbitrary $\tau\in\TT$ allows us to choose
\[
\tau=-F_n(\alpha)=-\frac{\alpha\rho_{n-1}(\alpha)}{\rho_{n-1}^*(\alpha)}
\Leftrightarrow
\tilde\tau=-\frac{\rho_n(\alpha)}{\rho_n^*(\alpha)}
\]
which fixes a node at $\alpha\in\TT$.

All this is well known but with our results of the previous section this can be generalized to the case where $S=\supp(\mu)$ is an arc on $\TT$.

\begin{theorem}[{\bf Szeg\H{o} and Szeg\H{o}-Radau q.f.}]\label{RadauSzego}
Let  $\mu$ be a positive measure supported on an arc $S=[a,b]\subset\TT$ and let $I_n(f)$ be a q.f.\ of the form (\ref{qf}).
\begin{enumerate}
\item $I_n(f)$ is a Szeg\H{o} quadrature formula exact in
    ${\Ll}_{2n-1}(\omega)$ if and only if its nodal polynomial is the $\tau$-invariant $(n,1)$-QPOPUC satisfying Theorem~\ref{zero} and $\omega=\tau\varphi_{-\delta_n}(\tau)$.
   \item $I_n(f)$ is the unique Szeg\H{o}-Radau quadrature
   formula exact in $\Ll_{2n-2}$ and prefixed node $\alpha\in\disk{S}\subseteq\TT$ if and only if its nodal polynomial is the $\tau$-invariant $(n,1)$-QPOPUC with $-\tau=\alpha\rho_{n-1}(\alpha)/\rho_{n-1}^*(\alpha)$ satisfying Theorem~\ref{zero}.
   If $\omega=\tau\varphi_{-\delta_n}(\tau)$, then the q.f.\ is exact in $\Ll_{2n-1}(\omega)$.
\end{enumerate}
\end{theorem}
\begin{proof}
The original proof of exactness in $\Ll_{2n-2}$ and the positivity of the weights that is given in \cite{JNT} for a measure supported on $\TT$ still holds if the support is an $\arc[a,b]\subset\TT$ and all the nodes are in the open $\arc(a,b)$.
We note that the positivity of the weights also follows directly from the exactness in $\Ll_{2n-2}$.
Indeed, as we mentioned above, by \cite[Theorem~2.3]{CDP}, the exactness in $\Ll_{2n-2}$ implies the positivity of $\lfloor\frac{2(n-1)}{2}\rfloor+1=n$ weights.
The statement about exactness in $\Ll_{2n-1}(\omega)$ is a consequence of
Theorem~\ref{Jacobi} and Lemma~\ref{lemomega}.
\end{proof}
\begin{remark}
The q.f.\ using $Q_{n,1}$ as nodal polynomial where we preselect $\tau$ and hence select some particular value of $\omega$ that defines the domain of exactness is the analogue of the Gaussian rule for a measure supported on $[a,b]\subset\RR$. If $\supp(\mu)=\TT$ this was considered in \cite{CDP}.

The q.f.\ using $Q_{n,1}$ as nodal polynomial where $\tau$ is fixed by the preselected node $\alpha$ is the analogue of the Gauss-Radau rule for a measure supported on $[a,b]\subset\RR$. Also this was considered in the case that $\supp(\mu)=\TT$ in e.g.\ \cite{Bul2}.

To the best of our knowledge our more general results are new when $\supp(\mu)$ is an arc $S\subsetneq \TT$.
\end{remark}

\subsection{Positive Szeg\H{o}-type q.f.\ with 2 or 3 prefixed nodes}

In the classical Gauss-Lobatto quadrature formulas fix zeros at the boundaries of $\supp(\mu)=[a,b]\subset\RR$.
For Szeg\H{o}-Lobatto q.f.\ when $\supp(\mu)=\TT$, there are no boundary points and therefore two distinct points $\alpha_1$ and $\alpha_2$ are preselected on $\TT$ and a positive q.f.\ having these nodes are constructed like in \cite{JaRei}.
For more details about Szeg\H{o}-Lobatto q.f.\ see also \cite{CDP,rog}.
In \cite{BGHN09} Szeg\H{o}-Lobatto quadrature is discussed in a more general context of orthogonal rational functions, but it contains the polynomials as a special case.

We summarize the results as follows.
\begin{theorem}[{\bf Szeg\H{o}-Lobatto quadrature formula}]\label{thmSL}
Let $n \geq 3$. Given $\tau\in\TT$ and distinct $\{\alpha_i\}_{i=1}^2\subset \TT$, we define $\{f_i\}_{i=1}^2\subset\TT$ as in (\ref{eqfi}). Assume the $f_i$ are distinct and $\alpha_1 f_1\ne \alpha_2 f_2$. Let $I_{n,3}$ be a q.f.\ of the form (\ref{qf}) with nodal polynomial $W_n$. Then $I_{n,3}$ is a Szeg\H{o}-Lobatto q.f exact in $\Ll_{2n-4}$
whose nodes include $\alpha_1$ and $\alpha_2$
if and only if
$W_n=Q_{n,3}$ a QPOPUC as in (\ref{eqQn3}) that satisfies the conditions of Theorem~\ref{pts2} or equivalently Theorem~\ref{pts2b}.
If moreover $\tau$ satisfies $\tau \overline \eta + \delta_{n-1} \neq 0$ with  $\eta$  given by (\ref{eta}), then $\omega=\frac{\tau\overline\eta+\delta_{n-1}}{\overline\tau\eta+\overline{\delta_{n-1}}}\in\TT$
is well defined and then $I_{n,3}$ is exact in $\Ll_{2n-1}(\omega)$.
\end{theorem}
\begin{proof}
It has been proved by
Theorems~\ref{pts2} and \ref{pts2b} under what conditions $Q_{n,3}$ has simple nodes on $\TT$ and by Theorem~\ref{Jacobi}(1) we know when
$I_{n,3}(f)$ is exact in $\Ll_{2n-4}$.
If $\tau\overline\eta+\delta_{n-1}\ne0$, then the orthogonality parameter $\tilde\tau$ and $\omega$ will be well defined by Lemma~\ref{lemomega} in the Appendix and the quadrature will be exact in $\Ll_{2n-3}(\omega)$ with $\omega=\tau\tilde\tau$ by Theorem~\ref{Jacobi}(2).

It only remains to prove that the weights are positive.
Therefore note that
\[
W_n(z) =  z\hat{\rho}_{n-1}(z)- \tau \hat{\rho}_{n-1}^*(z),~~\text{with}~~ \hat{\rho}_{n-1}(z)=z\rho_{n-2}(z)-\tau\overline{\eta}\rho_{n-2}^*(z),~~\tau\overline{\eta}\in\DD.
\]
By Favard's theorem on the unit circle (see \cite{JNT}), we conclude that $\hat{\rho}_{n-1}$ is a monic Szeg\H{o} polynomial and $W_n$ is a monic $(n,1)$-QPOPUC associated with a positive modified measure $\hat{\mu}$, then our $I_{n,3}(f)$ is also a Szeg\H{o}-Radau q.f.\ associated to $\hat{\mu}$ that coincides with our q.f.\ in $\mathcal{L}_{2n-3}(\omega)$, and hence it is a positive q.f.

So far we proved that the conditions of Theorems~\ref{pts2} or \ref{pts2b} imply the existence of the prescribed positive q.f.
The inverse implication is a consequence of \cite[Theorem 3.1]{Peh}.
\end{proof}

Since the previous Theorem~\ref{thmSL} holds for $\tau$ belonging to an arc, we may use it to fix a third node.

\begin{theorem} [\textbf{$(n,3)$-Szeg\H{o}-Peherstorfer quadrature}]\label{thmSP3}
Let $n\geq 3$ and let $I_{n,3}$ be the q.f.\ of the from (\ref{qf}) with nodal polynomial $W_n$.
Then $I_{n,3}$ is the unique Szeg\H{o}-Peherstorfer q.f.\ exact in $\Ll_{2n-4}$ with prescribed nodes $\alpha_1$, $\alpha_2$, and $\alpha_3$ if and only if
$W_n=Q_{n,3}$ as in (\ref{eqQn3}) satisfying the conditions of Theorem~\ref{3pointsfixed}.
If moreover $\tau$ verifying (\ref{autom}) satisfies $\tau \overline \eta + \delta_{n-1} \neq 0$ with  $\eta$  given by (\ref{eta}), then
$I_{n,3}$ is exact in $\Ll_{2n-3}(\omega)$ with $\omega=\frac{\tau\overline{\eta}+\delta_{n-1}}{\overline\tau\eta+\overline{\delta_{n-1}}}$.
\end{theorem}
\begin{proof}
By Theorem~\ref{3pointsfixed}, we know that if $\arc(\alpha_1,\alpha_2,\alpha_3)$ and $\arc(f_1,f_2,f_3)$ have the same orientation, then there exists a unique $Q_{n,3}$ with simple zeros on $\TT$
and $\alpha_1,\alpha_2$ and $\alpha_3$ are three of them.
Thus $I_{n,3}$ will be exact in ${\cal L}_{2n-4}$ by Theorem~\ref{Jacobi}(1). The weights are positive by the same argument as used in the previous theorem~\ref{thmSL}. The reciprocal is again a consequence of  \cite[Theorem 3.1]{Peh}.
The expression for $\omega$ follows from Theorem~\ref{Jacobi}  and Lemma~\ref{lemomega}.
\end{proof}

Theorems~\ref{thmSL} and \ref{thmSP3} hold again for integrals where $\supp(\mu)=\TT$.
If $\supp(\mu)$ is an arc $[a,b]\subsetneq\TT$, these theorems only guarantee $n-3$ zeros of odd multiplicity on the $\arc(a,b)$.
The `classical' versions of these Theorems do hold when two of the prefixed points are the boundary points $a$ and $b$.

\begin{theorem}[{\bf Classical Szeg\H{o}-Lobatto and Szeg\H{o}-Peherstorfer q.f.}]\label{CRadauSzego}

Let  $\mu$ be a measure supported on an arc $S=[a,b]\subsetneq\TT$.
Let $n\ge3$, $\tau \in \TT$, $f_a = F_{n-1}(a)$ and   $f_b = F_{n-1}(b)$  as in (\ref{eqfi})
and let $I_{n,3}$ be the q.f.\ of the form (\ref{qf}) with nodal polynomial $W_n$.
\begin{enumerate}
\item $I_{n,3}$ is a Szeg\H{o}-Lobatto quadrature formula (fixing nodes in $a$ and $b$) exact in
  ${\Ll}_{2n-3}(\omega)$ if and only if
  $W_n=Q_{n,3}$ as in (\ref{eqQn3})
  satisfying the first condition of Theorem~\ref{thm23}, with $\eta\in\DD$, $\tau \overline \eta + \delta_{n-1} \neq 0$
   and $\omega =\frac{\tau\overline{\eta}+\delta_{n-1}}{\overline\tau \eta+\overline{\delta_{n-1}}}$.
   Note that $\eta$ given by (\ref{eta}) is in $\DD$ if
  $\arc(a,\tau f_a,b)$ and $\arc(a,\tau f_a,\tau f_b)$ have opposite orientation.
\item $I_{n,3}$ is the unique Szeg\H{o}-Peherstorfer quadrature formula exact in $\Ll_{2n-4}$ and prefixed nodes $a$, $b$ and $\alpha\in\arc(a,b)$ if and only if $W_n=Q_{n,3}$ as in (\ref{eqQn3}) satisfying the second condition of Theorem~\ref{thm23} with $\eta\in\DD$, i.e.
    $\arc(a,\alpha, b)$ and $\arc(f_a,f_\alpha, f_b)$ have the same orientation, where $f_\alpha = \overline{F_{n-1}(\alpha)}$.
    If moreover  $\tau$  is such that $\tau\overline{\eta}+\delta_{n-1}\ne0$
then $I_{n,3}$ is exact in $\Ll_{2n-3}(\omega)$
with $\omega=\frac{\tau\overline{\eta}+\delta_{n-\ell}}{\overline\tau \eta+\overline{\delta_{n-\ell}}}$.
\end{enumerate}
\end{theorem}
\begin{proof}
By Theorems~\ref{thm23}, it is sufficient that the
nodal polynomial $Q_{n,3}$ satisfies the conditions mentioned, for both q.f.\ to have $n$ simple nodes in the arc $[a,b]\subset\TT$.
The positivity of the weights is deduced by the same argument used in the previous theorem.
The necessity of the conditions is a consequence of  \cite[Theorem 3.1]{Peh}.
\end{proof}

\subsection{Positive Szeg\H{o}-type q.f. with an arbitrary number of prefixed nodes}

For quadrature formulas with more than 3 prefixed nodes, we shall only consider integrals where $\supp(\mu)=\TT$ as we did in the previous section for the location of the zeros of QPOPUC of order $2\ell+1>3$.

\begin{theorem}[{\bf General Szeg\H{o}-Peherstorfer quadrature formula}]
\label{thmSPd}
Let $n\ge2\ell+1$, and
let ${I}_{n,2\ell+1}$ be a q.f.\ of the form (\ref{qf}) with nodal polynomial $W_n$.
\begin{enumerate}
\item
$I_{n,2\ell+1}$ is a Szeg\H{o}-Peherstorfer q.f.\  exact in
$\Ll_{2(n-\ell)-1}(\omega)$ with prefixed nodes
$\{\alpha_i\}_{i=1}^{2\ell}\subset\TT$
if and only if $W_n=Q_{n,2\ell+1}$ as in (\ref{Qn2lp1})
satisfying the hypothesis of Theorem~\ref{lemaadd} where $\tau$ is such that $\tau\overline{P_\ell(0)}-\delta_{n-\ell}\ne0$ and
$\omega=\frac{\tau\overline{P_\ell(0)}-\delta_{n-\ell}}{\overline\tau P_\ell(0)-\overline{\delta_{n-\ell}}}$.
\item
$I_{n,2\ell+1}$ is the unique Szeg\H{o}-Peherstorfer q.f.\ of the form (\ref{qf}) exact in
$\Ll_{2(n-\ell-1)}$ and prefixed nodes
$\{\alpha_i\}_{i=1}^{2\ell+1}\subset\TT$
if and only if
$W_n=Q_{n,2\ell+1}$ as in (\ref{Qn2lp1}) satisfying the hypothesis of Theorem~\ref{lemaadd2}. If moreover  $\tau$  is such that $\tau\overline{P_\ell(0)}-\delta_{n-\ell}\ne0$
then $I_{n,2\ell+1}$ is exact in $\Ll_{2(n-\ell-1)}(\omega)$
with $\omega=\frac{\tau\overline{P_\ell(0)}-\delta_{n-\ell}}{\overline\tau P_\ell(0)-\overline{\delta_{n-\ell}}}$.
\end{enumerate}
\end{theorem}
\begin{proof}
The existence and uniqueness of $W_n=Q_{n,2\ell+1}$ is the result of Theorem~\ref{lemaadd}
and Theorem~\ref{lemaadd2}, and the exactness in $\Ll_{2(n-\ell)-1}(\omega)$ is an immediate consequence of Theorem~\ref{Jacobi}.
As in the previous case, we are going to show that $Q_{n,2\ell+1}$ is a $(n,1)$-QPOPUC of a positive auxiliary measure.
Therefore, the associated q.f.\ is actually a Szeg\H{o}-Radau quadrature for this
auxiliary measure that coincides with our q.f.\ $I_{n,2\ell}(f)$ in $\Ll_{2(n-\ell)-1}(\omega)$, hence it is a positive q.f.
We know that the Schur-Cohn test for $P_\ell$ is equivalent with
\[
\left[\begin{array}{c} {P}_\ell(z)\\{P}_\ell^*(z)\end{array}\right] = \hat\theta_\ell \hat\theta_{\ell-1}~\cdots~\hat\theta_1
\left[\begin{array}{c} 1\\1\end{array}\right],~~\text{with}~~
\hat\theta_j=\left[\begin{array}{cc}1&\kappa_j\\\overline{\kappa_j}&1\end{array}\right]
\left[\begin{array}{cc}z&0\\0&1\end{array}\right],~~\kappa_j={P}_j(0)\in\DD,
\]
or, after multiplying the second row with $\tau$ and taking transpose
\[
[z~~\tilde\tau]
\theta'_1~\cdots~\theta'_\ell
\left[\begin{array}{cc}z&0\\0&1\end{array}\right]^{-1}=
[{P}_\ell(z)~~~\tilde\tau{P}_\ell^*(z)],~~
\theta'_j=
\left[\begin{array}{cc}1&\tilde\tau\overline{\kappa_j}\\\overline{\tilde\tau}{\kappa_j}&1\end{array}\right]
\left[\begin{array}{cc}z&0\\0&1\end{array}\right],
\]
while
\[
\left[\begin{array}{c} \rho_{n-\ell-1}(z)\\\rho_{n-\ell-1}^*(z)\end{array}\right] = \theta_{n-\ell-1}\theta_{n-\ell-2}~\cdots~\theta_1
\left[\begin{array}{c} 1\\1\end{array}\right], ~~\theta_j=\left[\begin{array}{cc}1&\delta_j\\\overline{\delta_j}&1\end{array}\right]
\left[\begin{array}{cc}z&0\\0&1\end{array}\right]
\]
where the $\delta_j$ is the $j$th Schur parameter.

Now define
\[
\left[\begin{array}{c}\tilde\rho_{n-\ell}(z)\\\tilde\rho^*_{n-\ell}(z)\end{array}\right]=
\tilde\theta_{n-\ell}
\left[\begin{array}{c}\rho_{n-\ell-1}(z)\\\rho_{n-\ell-1}^*(z)\end{array}\right],~~\text{with}~~
\tilde\theta_{n-\ell}=\theta'_\ell,
\]
and recursively
\[
\left[\begin{array}{c}\tilde\rho_{n-\ell+j}(z)\\\tilde\rho^*_{n-\ell+j}(z)\end{array}\right]=
\tilde\theta_{n-\ell+j}
\left[\begin{array}{c}\tilde\rho_{n-\ell+j-1}(z)\\\tilde\rho^*_{n-\ell+j-1}(z)\end{array}\right],~~j=1,\ldots,\ell-1
\]
where
\[
\tilde\theta_{n-\ell+j}=
\left[\begin{array}{cc}1&\tilde\delta_{n-\ell+j}\\\overline{\tilde\delta_{n-\ell+j}}&1\end{array}\right]=
\theta'_{d-j}=
\left[\begin{array}{cc}1&\tilde\tau\overline{\kappa_{\ell-j}}\\\overline{\tilde\tau}{\kappa_{\ell-j}}&1\end{array}\right].
\]
Thus, we end up with
\begin{eqnarray*}
Q_{n,d}(z)&=&[zP_\ell(z)~~\tilde\tau P_\ell^*(z)]\left[\begin{array}{c}\rho_{n-\ell-1}(z)\\\rho_{n-\ell-1}^*(z)\end{array}\right]\\
&=&
[z~~~\tilde\tau]\tilde\theta_{n-1}~\cdots~\tilde\theta_{n-\ell}\theta_{n-\ell-1}~\cdots~\theta_1 \left[\begin{array}{c}1\\1\end{array}\right]\\
&=&z\tilde\rho_{n-1}(z)+\tilde\tau\tilde\rho_{n-1}^*(z).
\end{eqnarray*}
In other words, the $\tilde\rho_{n-1}$ are monic Szeg\H{o} polynomials for a positive measure on the unit circle
by Favard's theorem since it has  parameters $\{\delta_1,\delta_2,\ldots,\delta_{n-\ell-1}\tilde\delta_{n-\ell},\ldots,\tilde\delta_{n-1}\}$ that are all in $\DD$.
Therefore all the weights are positive.
This shows the existence of a positive q.f.\ if all the conditions of Theorem~\ref{lemaadd} or Theorem~\ref{lemaadd2} are satified.
As in previous theorems, the reciprocal is a consequence of  \cite[Theorem~3.1]{Peh}.
\end{proof}

Note that F.~Peherstorfer in Theorem~3.1 of his paper \cite{Peh} characterized positive quadratures exact in the spaces $\mathcal{L}_{2(n-\ell-1)}$ in terms of quasi-paraorthogonal polynomials of odd order but without observing the orthogonality conditions that these verify.
However, he did not provide strategies for the construction of these positive q.f., and this has been precisely one of our purposes, taking advantage of the use of QPOPUC.
This is why we used the term Szeg\H{o}-Peherstorfer quadrature.

Our theorems recover these results but they are important because we also show
the influence of the invariance parameter on the nodes and on the maximal domain of exactness.
It opens the road to an analysis of the dynamics of the nodes as they vary in terms of this invariance parameter.
But more importantly,
our theorems are constructive and actually describe a computational
procedure to obtain the Schur coefficients needed to generate the nodal polynomial, and from these, it is possible to compute the nodes and weights of the  q.f.\ (see below in the numerical section).

If $2\ell+1$ nodes are fixed then there can be at most one solution for the $Q_{n,2\ell+1}$ and thus there is also a positive q.f.\ with these prefixed nodes and exact in $\Ll_{2(n-\ell-1)}$ will be unique if it exists.
It will also be exact in $\Ll_{2(n-\ell)-1}(\omega)$ with $\omega$ defined by the context.
If however we only fix $2\ell$ nodes, there may be an infinite set of
$Q_{n,2\ell+1}$, depending on its invariance parameter $\tau=Q_{n,2\ell+1}(0)$ that may vary in certain arcs on $\TT$.
The corresponding q.f., if it exists, is exact in
$\Ll_{2(n-\ell)-1}(\omega)$, where now $\omega$ depends on $\tau$.
The freedom to choose $\tau$ in certain arcs has been used before to fix an extra node.
We could use the same freedom to fix $\omega$ instead.
Given $2\ell+1$ nodes, find the corresponding
solution $Q_{n,2\ell+1}$, i.e., find $P_\ell$ and $\tau$ in (\ref{Qn2lp1}),
that matches these nodes gives a linear problem that we discussed above.
In essence $P_\ell$ is used to match $2\ell$ nodes and the parameter $\tau$ is used to match the last node.
Given $2\ell$ nodes and $\omega$, find the corresponding $Q_{n,2\ell+1}$ is the same when using $P_\ell$ to match the $2\ell$ nodes, but
finding a $\tau$ that gives a prescribed $\omega$ is a difficult nonlinear problem that does not allow an explicit solution.
However, generically there will be at most two values for $\tau$ as we prove in
Theorem~\ref{thm2om} of the Appendix.
In both cases, if a solution $Q_{n,2\ell+1}$ of these matching problems exists, then it remains to check that this results in a positive q.f. In practice, this is verified by the Schur-Cohn test for the polynomial $P_\ell$.

\section{Numerical experiments}\label{secNum}

The aim of this section is to present, some illustrative numerical examples in the construction of positive q.f.\ on the unit circle with a number of prescribed nodes greater than two and maximal domain of exactness. We consider the approximation of the integral (\ref{integral}) where $\mu$ is one of the most important measures on the unit circle: the Rogers-Szeg\H{o} (normalized) weight function (also known as ``wrapped" Gaussian measure), given by
\begin{equation}\label{rogers}
\mu(\theta,q)=\mu(\theta)=\frac{1}{\sqrt{2\pi\log{(1/q)}}}\sum_{j=-\infty}^\infty{\exp{\left(-\frac{(\theta-2\pi j )^2}{2\log{(1/q)}}\right)}} , \quad 0<q<1.
\end{equation}
Properties and applications of Rogers-Szeg\H{o} polynomials, the family of orthogonal polynomials on $\TT$ with respect to $\mu$ given by (\ref{rogers}), have been extensively studied in recent years.
In particular, quadrature rules associated with this measure were considered in \cite{rog}.

It is very well known that for this weight function, the trigonometric moments and monic Szeg\H{o} polynomials are explicitly given by $\mu_k=q^{\frac{k^2}{2}}$, $k \in \ZZ$, and by
\begin{equation}\nonumber 
\rho_n(z)=\sum_{j=0}^{n} (-1)^{n-j} \left[ \substack{n \\ j }
\right]_q q^{\frac{n-j}{2}}z^j,
\end{equation}
where for all $0< q < 1$, the usual $q$-binomial coefficients are defined by
\[
\left[ \substack{n \\ j } \right]_q=\frac{(n)_q}{(j)_q (n-j)_q}, \quad \textrm{where} \quad (n)_q=(1-q)(1-q^2)~\cdots~ (1-q^n) \quad \textrm{and} \quad (0)_q = \left[ \substack{n \\ 0 } \right]_q = \left[
\substack{n \\ n } \right]_q  \equiv 1,
\]
respectively. So, $\left[ \substack{n \\ j } \right]_q=\frac{(1-q^n)~\cdots~(1-q^{n-j+1})}{(1-q)~\cdots~(1-q^j)}$, $\delta_n=(-1)^{n}q^{\frac{n}{2}}$ and $\rho_n^*(z)=(-1)^n
q^{-\frac{n}{2}}\rho_n(qz)$.

The numerical computations were implemented in \textsf{Matlab} software with double precision.

Suppose first we want to construct a $16$-point rule with $6$ prescribed nodes ($d=2\ell=6$):
$\alpha_1=e^{-\frac{3\pi}{4}i}$, $\alpha_2= -i$, $\alpha_3=1 $, $\alpha_4=e^{\frac{\pi}{4}i}$, $\alpha_5=\overline{\alpha_2}$  and $\alpha_6=\overline{\alpha_1}$ and
we want to guarantee positive weights and a maximal degree of exactness, which is a space $\Ll_{25}(\omega)$ for some $\omega$.

We know $Q_{16,7}$ is given by
\[
Q_{16,7}(z) = zP_3(z)\rho_{12}(z)+
\tau_3P_3^*(z)\rho_{12}^{*}(z),
~~ P_3(z)=(z-\eta_1)(z-\eta_2)(z=\eta_3),
\]
and the maximal domain of exactness of the positive q.f.\ $I_{16,7}(f)$ is
\[
\Ll_{25}(\omega_{3})=\LL_{-12,12}\oplus \Span \left\{ z^{13} - \frac{\omega_3}{z^{13}}\right\},~~
\omega_3=\tau_3\tilde\tau_3,~~\tilde\tau_3=\frac{\overline{\eta_1\eta_2\eta_3}+\overline{\tau_3}\delta_{13}}
            {\overline{\tau_3}{\eta_1\eta_2\eta_3}+\overline{\delta_{13}}}.
\]
Observe that the nodes of the q.f.\ that were not prefixed, the $\{\eta_i\}_{i=1}^3$, the orthogonality parameter $\tilde\tau_3$ and thus also $\omega_3$ all depend on the choice for the invariance parameter $\tau_3$. We have chosen
$\tau_3=e^{0.9\pi i}$, resulting in $\tilde\tau\approx e^{-0.87834\pi i}$ and thus $\omega_3\approx e^{0.02166 \pi i}$
for the example in Table~\ref{nodesweights1} below.
The domain of exactness is maximized for this $\omega_3$.

\begin{figure}[!ht]
\begin{center}
\begin{tabular}{cc}
{\includegraphics[width=60mm]{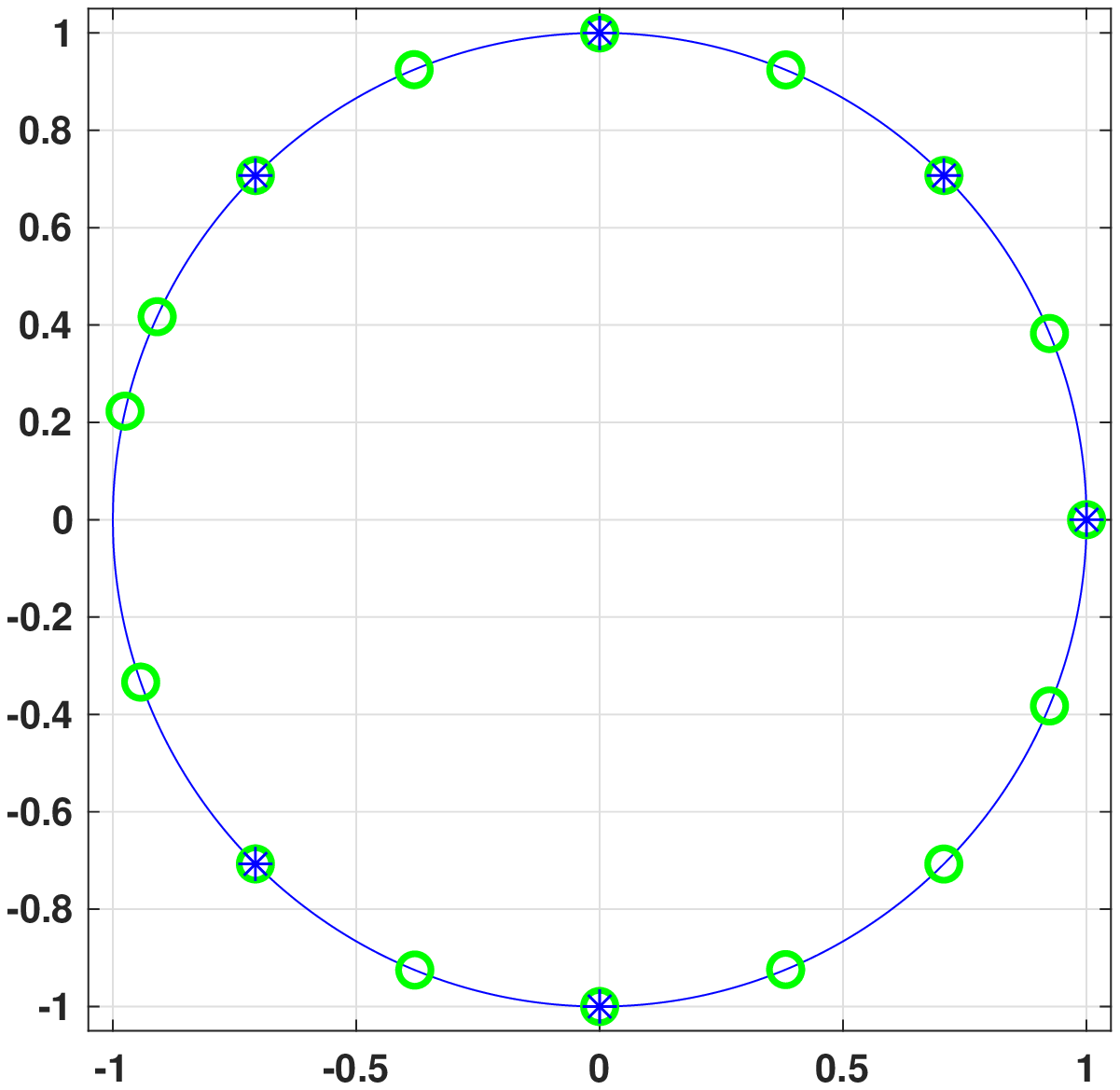}} &
{\includegraphics[width=60mm]{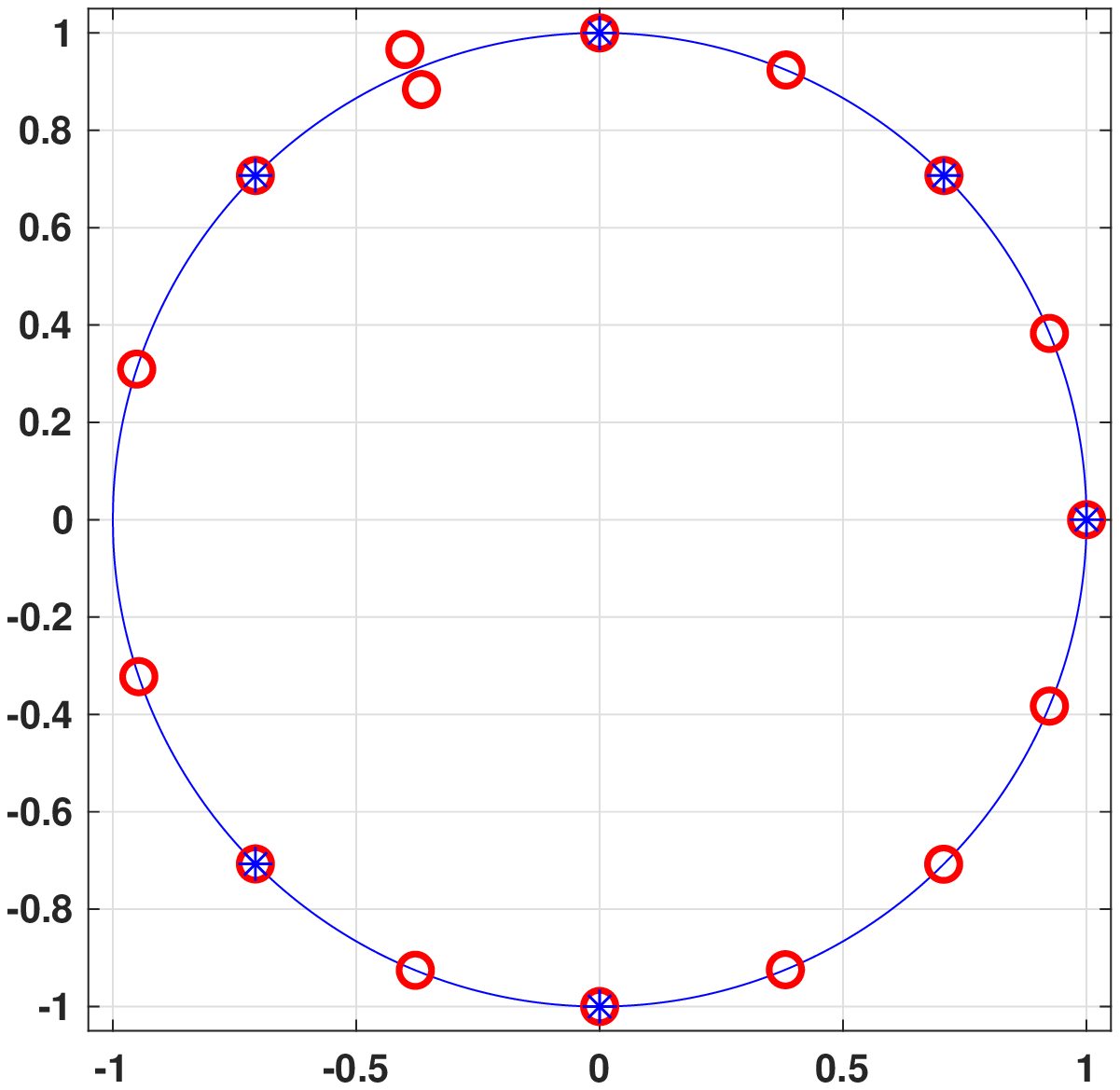}}\\
A & B\\
{\includegraphics[width=60mm]{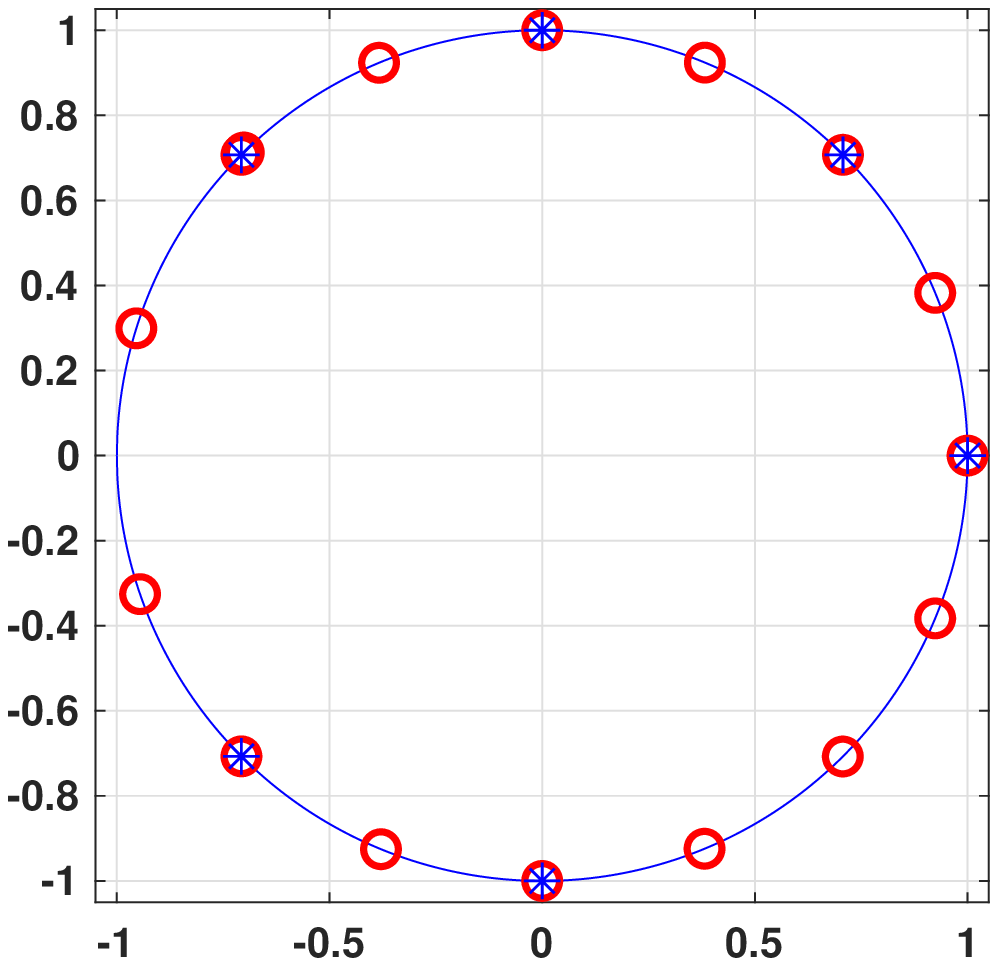}} &
{\includegraphics[width=60mm]{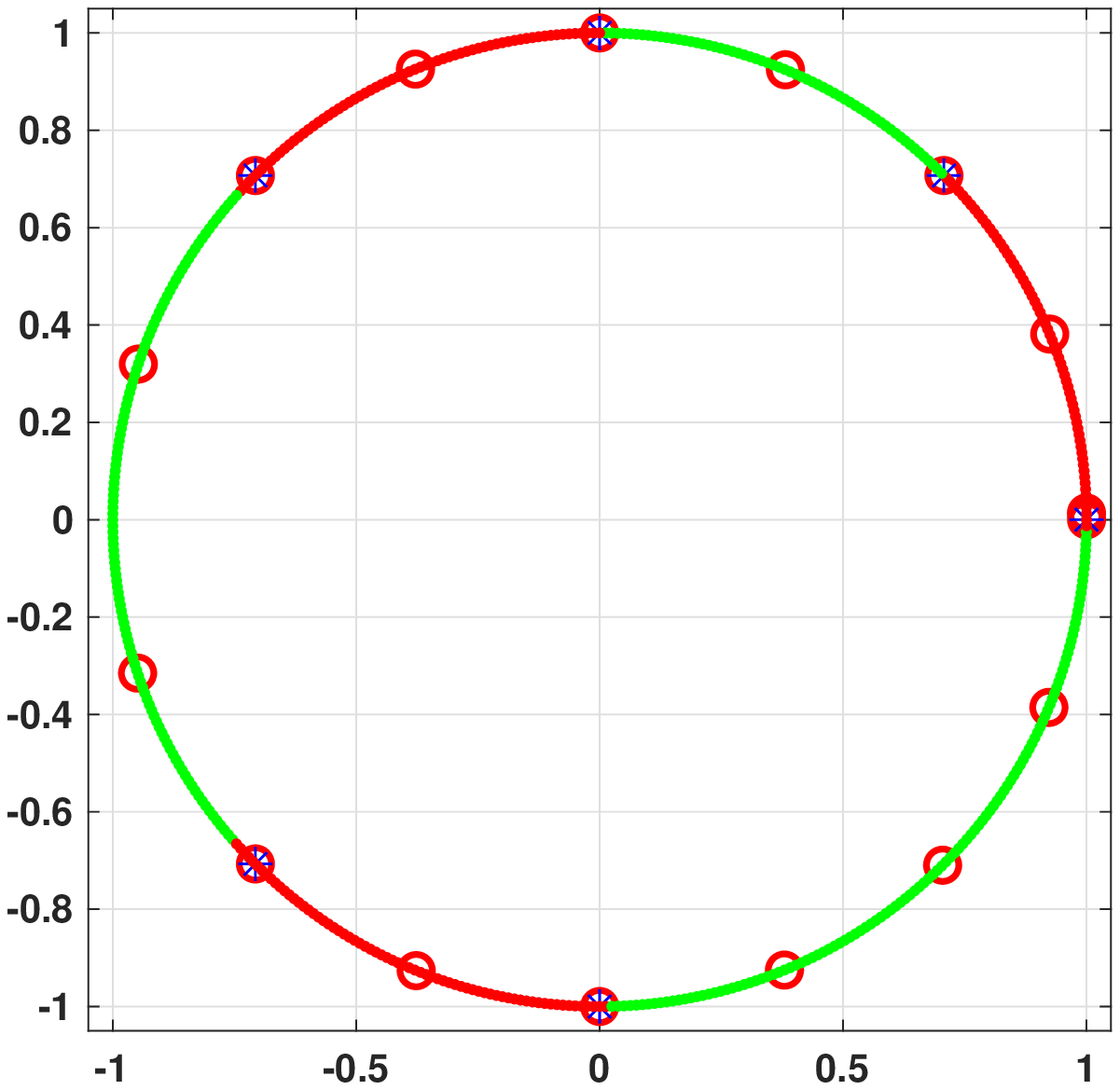}}\\
C & D
\end{tabular}
\vspace{.5cm} \caption{The plot shows the distribution of nodes in the Rogers-Szeg\H{o} example.
The circles show the locations of the computed nodes of the q.f.\ with the stars indicating the preselected nodes.
The figure corresponds to the nodal polynomial $Q_{16,7}$ with 6 preselected nodes.
The q.f.\ is exact in $\Ll_{25}(\omega)$.
(A): $\tau=\exp({0.9\pi i})$ corresponding to $\omega\approx \exp({0.02166 \pi i})$.
(B): $\tau=\exp({0.63\pi i})$ there is one pair of nodes not on $\TT$.
(C): $\tau\approx \exp({0.7575\pi i})$ one moving zero coincides with one of the prefixed (in this case with $\exp({0.75\pi i})$).
(D): Nodes for $\tau=1$. There are 16 simple nodes, but not all the weights are positive. For $\tau$ in a green arc then the weights are positive. For $\tau$ in a red arc the weights are not positive.
} \label{fig1}
\end{center}
\end{figure}

Using the procedure developed in Theorem~\ref{thmSPd} we can compute the Schur parameters
$\tilde\delta_{13},\tilde\delta_{14},\tilde\delta_{15}$ defined by the preselected nodes.
These parameters, together with the Schur parameters $\delta_j$, $j=1,\ldots,12$, for the Rogers-Szeg\H{o} polynomials
allow us to compute the nodes and weights from the CMV matrix associated with them.
See \cite{CMV02,CMV06,CMP,ArtBCCB14}.
The results for $\tau_3=e^{0.9\pi i}$ are displayed in Table~\ref{nodesweights1} and the nodes are also graphically represented
in Figure~\ref{fig1}(A).

\begin{table}[!ht]
\renewcommand{\arraystretch}{1.4}
\[
\begin{array}{|c|c|}
\hline
\textrm{Nodes}&\textrm{Weights}\\\hline
  -0.942694568084626 - 0.333656936543722i & 0.000883914413545\\\hline
 \alpha_1=\frac{\sqrt{2}}{2}(-1-i)       &  0.003573449079563\\\hline
 -0.379700471214962 - 0.925109481174597i & 0.011670163240250\\\hline
 \alpha_2=-i                             & 0.031571926071034\\\hline
  0.382004048856982 - 0.924160649809800i & 0.069097384838415\\\hline
  0.706940179145082 - 0.707273343984008i & 0.120722316200665\\\hline
  0.923913888542018 - 0.382600479036773i & 0.168408825865268\\\hline
  \alpha_3=1                             & 0.188077141674534\\\hline
  0.923939153768060 + 0.382539462192281i & 0.168360171656870\\\hline
  \alpha_4=\frac{\sqrt{2}}{2}(1+i)	 & 0.120696119860582\\\hline
  0.382340759272763 + 0.924021397911716i & 0.069140956573729\\\hline
 \alpha_5= i				 & 0.031642200537795\\\hline
 -0.381050525980387 + 0.924554215095074i & 0.011673845455142\\\hline
 \alpha_6=\frac{\sqrt{2}}{2}(-1+i)	 & 0.003486160898429\\\hline
 -0.908886480405499 + 0.417043601720617i & 0.000537635619058\\\hline
 -0.974806799005855 + 0.223050901392394i & 0.000457788015119\\\hline
 \end{array}
 \]
\caption{A $16$-point q.f.\ for the Rogers-Szeg\H{o} weight function with prescribed nodes $\{ \alpha_i \}_{i=1}^{6}$ and maximal domain of exactness $\Ll_{25}(\omega)$.
We have taken $q=0.5$ and $\tilde\tau_3=e^{0.9\pi i}$.}\label{nodesweights1}
\renewcommand{\arraystretch}{1.5}
\end{table}

As we have mentioned before, positive quadrature formulas will be guaranteed for $\tau_3$ in certain arcs of $\TT$.
In Figure~\ref{fig1}(D) we have indicated these as green arcs.
In our example, these are the arcs given approximately by
$\{e^{\theta\pi i}: \theta\in(0.251,0.499)\cup(0.765,1.229)\cup(1.505,1.995)\}$.
For each $\tau$ in the three green arcs, the conditions of Theorems~\ref{lemaadd} and \ref{thmSPd} are satisfied
and we can obtain a positive q.f. with the 6 nodes prefixed
and we obtain exactness in $\Ll_{23}(\omega_3)$.

Note that the boundaries of the arcs seem to coincide with the fixed nodes,
but that is not true.
We know that at such a boundary point, one of the $\eta$'s crosses $\TT$
and at the moment it is on $\TT$, it must coincide with (at least one) zero of $Q_{16,7}$. Experimentally we observed that this $\eta$ leaves $\TT$ by passing through one of the prefixed $\alpha$'s. This creates some `degenerate situation' where one of the `moving zeros' of $Q_{16,7}$, (i.e.\ the zeros varying with $\tau$) collides with that $\alpha$ (as in Figure~\ref{fig1}(C) where the collision happens in $\alpha=\exp(0.75\pi i)$ and only 15 distinct nodes are shown).
Since that $\alpha$ is then a double zero, it will also be a zero of the derivative $Q'_{16,7}$, while all the other zeros of $Q'_{16,7}$ stay in $\DD$.
It is intuitively clear that the `moving zeros' closest to $\tau$ will be most sensitive to varying $\tau$. That explains why the collision of a moving zero and a prefixed $\alpha$ will happen near the value of $\tau$, and thus why the boundary points for the arcs of $\tau$ will be near the prefixed zeros.

The complementary red arcs correspond to situations where at least one of the $\eta$'s, and hence one of the Schur parameters for $P_3$ is not in $\DD$.
When $\tau$ enters a red arc, no positive quadrature exist for several reasons.
What happens there depends on the location of the zeros of the derivative $Q'_{16,7}$.
If all the zeros of $Q'_{16,7}$ are in $\DD$, then the zeros of $Q_{16,7}$ are simple and on $\TT$.
If $Q'_{16,7}$ has a zero in $\EE$, then there is at least one
pair of zeros of $Q_{16,7}$ not on $\TT$ as in
Figure~\ref{fig1}(B).
Since a zero of $Q'_{16,7}$ can only leave $\DD$ by crossing $\TT$, it can only pass $\TT$ through a multiple zero of $Q_{16,7}$. We observe that this is in a red arc where two moving zeros collide.

The remaining Figure~\ref{fig1}(D) shows the green and red arcs for $\tau$.
The nodes are plotted for the situation $\tau=1$ where there are 16 simple nodes on $\TT$ but not all the weights are positive.
As we have pointed out in Remark~\ref{remtaul}, for a deeper exploration it is necessary to make a detailed study of the behavior of the corresponding Schur parameters.

\begin{figure}[!ht]
\begin{center}
\begin{tabular}{c}
\includegraphics[width=8cm]{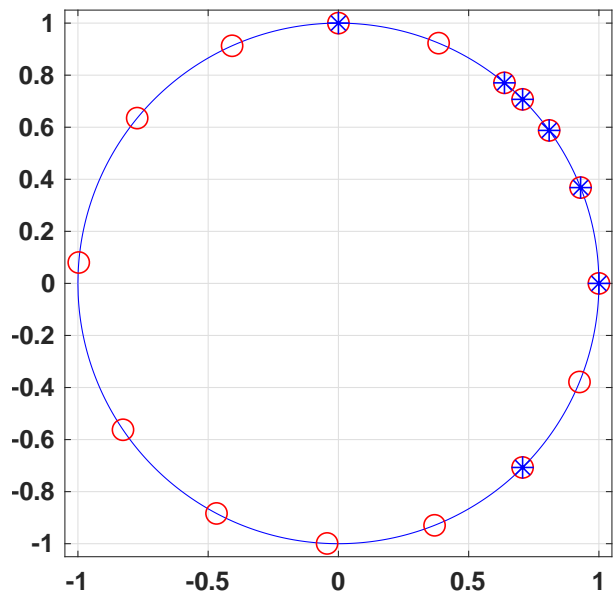}
\end{tabular}
\vspace{.5cm} \caption{The plot shows the distribution of nodes in the Rogers-Szeg\H{o} example.
Seven nodes are prefixed based on the nodal polynomial $Q_{16,7}$ with an appropriate choice of $\tau_3$
The red circles show the locations of the computed nodes of the q.f.\ with the stars indicating the preselected nodes.
The quadrature is exact in $\Ll_{25}(\omega_3)$ with $\omega_3\approx e^{0.8\pi i}$.
} \label{fig2}
\end{center}
\end{figure}

The freedom to choose the $\tau$ in these green arcs allows us to fix a seventh node.
This is described by Theorems~\ref{lemaadd2} and \ref{thmSPd}.
As an example we choose again the Rogers-Szeg\H{o} polynomials but now select seven nodes
that are less symmetric on the circle.
The results are shown for the nodal polynomial $Q_{16,7}$ in Table~\ref{tab2} and Figure~\ref{fig2}.

\begin{table}[!ht]
\renewcommand{\arraystretch}{1.4}
\[
\begin{array}{|r|c|}
\hline
\textrm{Nodes}&\textrm{Weights}\\\hline
 -0.826939158283098 - 0.562291408878033i &  0.000883914413545\\\hline
 -0.467939713417261 - 0.883760388684045i &  0.003573449079563\\\hline
 -0.043436557033519 - 0.999056187365392i &  0.011670163240250\\\hline
  0.369113471779490 - 0.929384336510408i &  0.031571926071034\\\hline
  \color{blue}{*~~0.707106781186549 - 0.707106781186549i} &  \color{blue}{0.069097384838415}\\\hline
  0.925520203512843 - 0.378698234600514i &  0.168718787427850\\\hline
  \color{blue}{*~~1.000000000000000 - 0.000000000000001i} &  \color{blue}{0.184031304001781}\\\hline
  \color{blue}{*~~0.929776485888252 + 0.368124552684676i} &  \color{blue}{0.157535679739229}\\\hline
  \color{blue}{*~~0.809016994374943 + 0.587785252292476i} &  \color{blue}{0.029784407222365}\\\hline
  \color{blue}{*~~0.707106781186542 + 0.707106781186551i} &  \color{blue}{0.087805627362136}\\\hline
  \color{blue}{*~~0.637423989748693 + 0.770513242775784i} &  \color{blue}{0.017212913090585}\\\hline
  0.384703682295919 + 0.923040127420233i &  0.067644558942418\\\hline
  \color{blue}{*~~0.000000000000000 + 1.000000000000000i} &  \color{blue}{0.032843677513034}\\\hline
 -0.408105016619360 + 0.912934989695385i &  0.011987419059172\\\hline
 -0.772470256507957 + 0.635050945051283i &  0.003103537707867\\\hline
 -0.996767928597351 + 0.080334902251429i &  0.000698363524723\\\hline
 \end{array}
 \]
\caption{A $16$-point q.f.\ for the Rogers-Szeg\H{o} weight function with prescribed nodes $\{ \alpha_i \}_{i=1}^{7}$ (in blue and marked by *) and maximal domain of exactness.
The required $\tau_3$ was $\tau_3=-0.363884303133021 - 0.931444155026694i\approx e^{-0.6\pi i}$ and the formula is exact in
$\mathcal{L}_{25}(\omega_3)$ with $\omega_3\approx e^{-0.8\pi i}$.
}\label{tab2}
\renewcommand{\arraystretch}{1.5}
\end{table}

\clearpage

\section*{Appendix}
\setcounter{section}{1}
\setcounter{theorem}{0}
\renewcommand{\thesection}{\Alph{section}}

In this addendum we collect some auxiliary properties that are used to prove some results in the body of this paper.

The first Lemma relates to points on $\TT$ and it is quoted from \cite[Chapter~3]{Ahlfors}.
\begin{lemma}\label{signo}
Let $[z_1,z_3]$ be an arc on $\TT$ and $z_2\in(z_1,z_3)$. Then
$$\frac{(z-z_1)(z_2-z_3)}{(z-z_3)(z_2-z_1)} \in \RR \Leftrightarrow z \in \TT.$$
Moreover,
$$\frac{(z-z_1)(z_2-z_3)}{(z-z_3)(z_2-z_1)} \in \RR^+ \Leftrightarrow z \in (z_1,z_3),$$
\end{lemma}

The next lemma is a simple geometric observations.

\begin{lemma}\label{chords}
Let $\alpha,\beta,\gamma,\delta$ be four distinct points on $\TT$ with $\alpha\delta-\beta\gamma\ne0$ then the secants $\alpha\beta$ and $\gamma\delta$ intersect in $\DD$ if and only if $\arc(\alpha,\beta,\gamma)$ and $\arc(\alpha,\beta,\delta)$ have opposite orientation.
\end{lemma}
\begin{proof}
The condition $\alpha\delta-\beta\gamma\ne0$ implies that the two secants are not parallel.
The secant $\alpha\beta$ divides $\TT$ in two arcs $A_1$ and $A_2$. It is obvious that the secant $\gamma\delta$ will intersect the secant $\alpha\beta$ in $\DD$ if and only if both $\gamma$ and $\delta$ do not belong to the same arc which may be either $A_1$ or $A_2$. This is equivalent with saying that $\arc(\alpha,\beta,\gamma)$ and $\arc(\alpha,\beta,\delta)$ have opposite orientation.
Similarly, it holds that the secants intersect in $\EE$ if and only if  $\arc(\alpha,\beta,\gamma)$ and $\arc(\alpha,\beta,\delta)$ have the same orientation.
\end{proof}

The following relations are easy to verify.
\begin{lemma}\label{lemA}
Let $n\ge 2\ell+1$.
Let
$Q_{n,2\ell+1}(z)=zp_\ell(z)\rho_{n-\ell-1}(z)+\tau p_\ell^*(z)\rho_{n-\ell-1}^*(z)$ with $p_\ell\in\PP_\ell\backslash\PP_{\ell-1}$ monic and $\tau\in\TT$.\\
Define $a_k=\langle Q_{n,2\ell+1},z^{k}\rangle$, $k=0,\ldots,n$, then
$a_k=\tau\overline{a_{n-k}}$ and $a_k=0$, $k=\ell+1,\ldots,n-\ell-1$ while
\[
a_{n-\ell}=({p_\ell(0)}-\tau\overline{\delta_{n-\ell}})\|\rho_{n-\ell-1}\|^2
\text{~~and~~}
\langle Q_{n,2\ell+1},z^{n-\ell}-\tilde\tau z^\ell\rangle
=
a_{n-\ell}-\overline{\tilde\tau}\tau \overline{a_{n-\ell}}.
\]
If $a_{n-\ell}=0$ or equivalently $\sigma_{n,\ell}=\overline{p_\ell(0)}-\overline\tau{\delta_{n-\ell}}=0$ then
there exists a monic $p_{\ell-1}\in\PP_{\ell-1}\backslash\PP_{\ell-2}$ such that
\[
Q_{n,2\ell+1}(z)=zp_{\ell-1}\rho_{n-\ell}(z)+\tau p_{\ell-1}^*(z)\rho_{n-\ell}^*(z)=Q_{n,2\ell-1}(z).
\]
If $a_{n-\ell}\ne0$ then
\[
\langle Q_{n,2\ell+1},z^{n-\ell}-\tilde\tau z^\ell\rangle=0 \Leftrightarrow
\tilde\tau=\tau\frac{\overline{a_{n-\ell}}}{a_{n-\ell}}=
\frac{\tau\overline{p_\ell(0)}-\delta_{n-\ell}}{{p_\ell(0)}-\tau\overline{\delta_{n-\ell}}}.
\]
Let $\tilde{Q}_{n,2\ell+1}(z)=q_\ell^*(z)\rho_{n-\ell}(z)+\tilde\tau q_\ell(z)\rho_{n-\ell}^*(z)$, with $q_\ell\in\PP_{\ell}\backslash\PP_{\ell-1}$ monic and $\tilde\tau\in\TT$.\\
Define
$\tilde{a}_k=\langle \tilde{Q}_{n,2\ell+1}(z),z^{k}\rangle$ then $\tilde{a}_k=\tilde\tau\overline{\tilde{a}_{n-k}}$, $k=0,\ldots,n$, and $\tilde{a}_k=0$ for $k=\ell+1,\ldots,n-\ell-1$ while
$\tilde{a}_{n-\ell}=\|\rho_{n-\ell}\|^2\ne0$, thus
$\tilde{Q}_{n,2\ell+1}\in\QQ_{n,2\ell+1}^\perp\backslash\QQ_{n,2\ell-1}^\perp$.\\
If $\nu_{n,\ell}=\overline{q_\ell(0)}+\tilde\tau\overline{\delta_{n-\ell}}=0$
then $\tilde{Q}_{n,2\ell+1}\in \PP_{n-1}$ and $\tilde{Q}_{n,2\ell+1}(0)=0$ and
thus $\tilde{Q}_{n,2\ell+1}(z)=zQ_{n-2,2\ell-1}(z)$ with $\langle Q_{n-2,2\ell-1},z^{n-\ell-1}\rangle\ne0$.\\
If  $\nu_{n,\ell}\ne0$ then
$\tilde{Q}_{n,2\ell+1}\in\QQ_{n,2\ell+1}^\perp\backslash\QQ_{n,2\ell-1}^\perp$
is in $\PP_n\backslash\PP_{n-1}$,
and if $\sigma_{n,\ell}\ne0$ then the monic polynomial
${Q}_{n,2\ell+1}\in\QQ_{n,2\ell+1}^\perp\backslash\QQ_{n,2\ell-1}^\perp$.
\\
Thus if $\nu_{n,\ell}\sigma_{n,\ell}\ne0$ then
$\nu_{n,\ell}^{-1}\tilde{Q}_{n,2\ell+1}(z)=Q_{n,2\ell+1}$ are equivalent representations of all monic polynomials in $\QQ_{n,2\ell+1}^\perp\backslash\QQ_{n,2\ell-1}^\perp$.
\end{lemma}
Note that $\tilde{a}_{n-\ell}=\tilde\tau\overline{\tilde{a}_\ell}>0$
so that $\langle \tilde{Q}_{n,2\ell+1},z^{n-\ell}-\tilde\tau z^\ell\rangle =0$.

The next result gives the relations between $\tau,\tilde\tau, p_\ell$ and $q_\ell$ for monic polynomials in $\QQ_{n,2\ell+1}^\perp\backslash\QQ_{n,2\ell-1}^\perp$.
\begin{lemma}\label{lemomega}
Let $p_\ell$ and $q_\ell$ be monic polynomials of degree $n$ and
assume $\nu_{n,\ell}=\overline{q_\ell(0)}+\tilde\tau\overline{\delta_{n-\ell}}\ne0$,
and
$\sigma_{n,\ell}=\overline{p_\ell(0)}-\overline\tau{\delta_{n-\ell}}\ne0$
then
\begin{eqnarray*}
Q_{n,2\ell+1}(z)&=& \nu_{n,\ell}^{-1}[q_\ell^*(z)\rho_{n-\ell}(z)+\tilde\tau q_\ell(z)\rho_{n-\ell}^*(z)],\\
&=&zp_\ell(z)\rho_{n-\ell-1}(z)+\tau p_\ell^*(z)\rho_{n-\ell-1}^*(z),
\end{eqnarray*}
represent all monic $(n,2\ell+1)$-QPOPUC not in $\QQ_{n,2\ell-1}^\perp$ and
\begin{eqnarray*}
&&p_\ell=\nu_{n,\ell}^{-1}[q_\ell^*+\tilde\tau \overline{\delta_{n-1}}q_\ell],~~
\nu_{n,\ell}=\overline{q_\ell(0)} +\tilde\tau \overline{\delta_{n-\ell}}\\
&&q_\ell=\sigma_{n,\ell}^{-1}[p_\ell^*-\overline\tau\delta_{n-\ell}p_\ell],~~
\sigma_{n,\ell}=\overline{p_\ell(0)}-\overline\tau{\delta_{n-\ell}}\\
&&{\nu_{n,\ell}}\overline{\sigma_{n,\ell}}={\overline{\nu_{n,\ell}}}{\sigma_{n,\ell}}=1-|\delta_{n-\ell}|^2>0\\
&&\tau=\tilde\tau\overline{\nu_{n,\ell}}/\nu_{n,\ell}=\tilde\tau\overline{\sigma_{n,\ell}}/{\sigma_{n,\ell}}\\
&& \omega=\tau\tilde\tau=\frac{ \tilde\tau q_\ell(0)+\delta_{n-\ell}}{\overline{\tilde\tau}\overline{q_\ell(0)}+\overline{\delta_{n-\ell}}}=
\frac{\tau\overline{p_\ell(0)}-\delta_{n-\ell}}{\overline{\tau}{p_\ell(0)}-\overline{\delta_{n-\ell}}}.
\end{eqnarray*}
\end{lemma}
\begin{proof}
Using the Szeg\H{o} recursion we get
\begin{eqnarray*}
Q_{n,2\ell+1}&=&\nu_{n,\ell}^{-1}[z \rho_{n-\ell-1}+\delta_{n-\ell}\rho_{n-\ell-1}^*]q_\ell^*+\tilde\tau
\nu_{n,\ell}^{-1}[\rho_{n-\ell-1}^*+\overline{\delta_{n-\ell}}z\rho_{n-\ell-1}]q_\ell\\
&=&\nu_{n,\ell}^{-1}[q_\ell^*+\tilde\tau\overline{\delta_{n-\ell}}q_\ell]z\rho_{n-\ell-1}+
\nu_{n,\ell}^{-1}[q_\ell^*\delta_{n-\ell}+\tilde\tau q_\ell]\rho_{n-\ell-1}^*\\
&=& zp_\ell\rho_{n-\ell-1}+\tau p_\ell^*\rho_{n-\ell-1}^*,~~p_\ell=\nu_{n,\ell}^{-1}[q_\ell^*+\tilde\tau\overline{\delta_{n-\ell}}q_\ell]
\end{eqnarray*}
with $\nu_{n,\ell}=\overline{q_\ell(0)} +\tilde\tau \overline{\delta_{n-\ell}}$
because $q_\ell$ is monic
and $\tau=\tilde\tau\overline{\nu_{n,\ell}}/\nu_{n,\ell}$.
For the inverse relation we use
\begin{eqnarray*}
p_\ell^*-\overline\tau\delta_{n-\ell}p_\ell &=&
\overline{\nu_{n,\ell}^{-1}}(q_\ell+\overline{\tilde\tau}\delta_{n-\ell}q_\ell^*)-\overline\tau \delta_{n-\ell}\nu_{n,\ell}^{-1}(q_\ell^*+\tilde\tau\overline{\delta_{n-\ell}}q_\ell)\\
&=& (\overline{\nu_{n,\ell}^{-1}}-\overline\tau\delta_{n-\ell}\nu_{n,\ell}^{-1}\tilde\tau\overline{\delta_{n-\ell}})q_\ell+
(\overline{\nu_{n,\ell}^{-1}}\overline{\tilde\tau}\delta_{n-\ell}-\overline\tau\delta_{n-\ell}\nu_{n,\ell}^{-1})q_\ell^*\\
&=&
\overline{\nu_{n,\ell}^{-1}}\left(1-\overline\tau\tilde\tau\frac{\nu_{n,\ell}^{-1}}{\overline{\nu_{n,\ell}^{-1}}}|\delta_{n-\ell}|^2\right)q_\ell+
\delta_{n-\ell}(\overline{\nu_{n,\ell}^{-1}}\overline{\tilde\tau}-\overline\tau \nu_{n,\ell}^{-1})q_\ell^*\\
&=&
\overline{\nu_{n,\ell}^{-1}}(1-|\delta_{n-\ell}|^2)q_\ell,
\end{eqnarray*}
where the last line follows from $\tau=\tilde\tau\frac{\nu_{n,\ell}^{-1}}{\overline{\nu_{n,\ell}^{-1}}}$.
Because $q_\ell$ is monic, ${\sigma_{n,\ell}}=\overline{\nu_{n,\ell}^{-1}}(1-|\delta_{n-\ell}|^2)=\overline{p_\ell(0)}-\overline\tau{\delta_{n-\ell}}$.
The remaining equalities are immediately verified.
\end{proof}

Let $I_{n,2\ell+1}$ be a positive Szeg\H{o}-Peherstorfer quadrature formula as described in Theorem~\ref{lemaadd2} with $2\ell$ prefixed nodes and exact in $\Ll_{2(n-\ell)+1}(\omega)$ with $\omega$ depending on $\tau$ as described in the theorem. Given $\tau$, we can compute $P_\ell$ and $\omega$.
The converse of finding $\tau$, given $\omega$ is a highly nonlinear problem.
However, generically, there will be two solutions as follows by the next theorem.
\begin{theorem}\label{thm2om}
Let $P_\ell$ be defined by the system (\ref{system2l}) parametrized in $\tau$ then for $\omega\in\TT$, the nonlinear equation
$\omega=\frac{\tau\overline{P_\ell(0)}-\delta_{n-\ell}}{\overline\tau P_\ell(0)-\overline{\delta_{n-\ell}}}$, will in general have at most two different solutions for $\tau\in\TT$.
\end{theorem}
\begin{proof}
The system (\ref{system2l}) can be written as
\[
\mathcal{M}\left[\begin{array}{cc}I&\\&\tau I \end{array}\right] \left[\begin{array}{c}\v{p}\\\overline{\v{p}}\end{array}\right]=
-\tau\v{f}-\v{d}.
\]
If $\det\mathcal{M}\ne0$, then
\[
\left[\begin{array}{c}\v{p}\\\overline{\v{p}}\end{array}\right]=
-\left[\begin{array}{cc}\tau I&\\& I \end{array}\right]
\mathcal{M}^{-1}\v{f}-
\left[\begin{array}{cc}I&\\&\overline\tau I \end{array}\right] \mathcal{M}^{-1}
\v{d}.
\]
The element on the first row is ${P_\ell(0)}=\prod_{i=1}^\ell(-{\eta_i})$ which has the form $A\tau+B$ with $[A,B]=-[1~0~\cdots~0]\mathcal{M}^{-1}[\v{f},\v{d}]$.
Substituting this into the equation for $\omega$ and writing $D$ for $\overline{\delta_{n-\ell}}$ gives
$\omega[\overline\tau(A\tau+B)-{D}]=\tau(\overline{A}\overline{\tau}+\overline{B})-\overline{D}$.
Since we are looking for $\tau\in\TT$, this results in a constrained quadratic equation
\[
\overline{B}\tau^2+\tau[(\overline{A}-\overline{D})-\omega(A-D)]-\omega B = 0,\quad \tau\in\TT.
\]
Thus except when $B=0$, which means that $\omega=\frac{\overline{A}-\overline{D}}{A-D}$ does not depend on $\tau$, and thus any $\tau$ will do,
there will be at most two solutions for $\tau\in\TT$.
\end{proof}

\bibliography{QPA}
\bibliographystyle{plain}

\end{document}